\documentclass[11pt]{amsart}
\usepackage{amsfonts}
\usepackage{amscd}
\usepackage{amssymb}
\usepackage{graphics}
\usepackage{amsmath}
\usepackage[english]{babel}
\usepackage{hyperref}
\usepackage{graphicx}%
\hypersetup{colorlinks,citecolor=blue}
\input xy
\xyoption{all}

\newcommand{\BZ}{{\mathbb{Z}}}

\newcommand{\BR}{{\mathbb{R}}}
\newcommand{\BC}{{\mathbb{C}}}
\newcommand{\BF}{{\mathbb{F}}}
\newcommand{\BQ}{{\mathbb{Q}}}
\newcommand{\BP}{{\mathbb{P}}}

\newcommand{\BG}{{\mathbb{G}}}

\newcommand{\BK}{{\mathbb{K}}}

\newcommand{\OO}{{\mathcal{O}}}
\newcommand{\gD}{\Delta}
\newcommand{\gd}{\delta}
\newcommand{\gb}{\beta}

\newcommand{\gC}{\Gamma}
\newcommand{\gc}{\gamma}
\newcommand{\gs}{\sigma}
\newcommand{\gS}{\Sigma}

\newcommand{\gl}{\lambda}\newcommand{\gL}{\Lambda}
\newcommand{\ga}{\alpha}

\def\bbr{{\Bbb R}}

\newcommand{\ti}[1]{\tilde{#1}}

\newcommand{\rank}{\text{rank}}
\newcommand{\Ad}{\text{Ad}}

\newcommand{\SL}{\text{SL}}
\newcommand{\GL}{\text{GL}}
\newcommand{\PSL}{\text{PSL}}
\newcommand{\PGL}{\text{PGL}}

\def\Aut{\text{Aut}}

\newtheorem{prop}{Proposition}[section]
\newtheorem{thm}[prop]{Theorem}
\newtheorem{lem}[prop]{Lemma}
\newtheorem{cor}[prop]{Corollary}
\newtheorem{conj}[prop]{Conjecture}
\theoremstyle{definition}

\newtheorem{defn}[prop]{Definition}
\newtheorem{rem}[prop]{Remark}

\newtheorem{exam}[prop]{Example}

\newtheorem{clm}[prop]{Claim}
\begin{document}

\author{Michelle Bucher and Tsachik Gelander}

\thanks{Michelle Bucher acknowledges support from the Swedish Research Council (VR) grant 621-2007-6250.
Tsachik Gelander acknowledges the financial support from the European Community's seventh Framework Programme (FP7/2007-2013) / ERC grant agreement 260508, and a partial support from the Israeli Science Foundation and the Gustafsson Foundation.}

\date{\today}
\title[The generalized Chern conjecture for $\mathcal{H}^n$-manifolds]{The generalized Chern conjecture for manifolds that are locally a product of surfaces} 
\begin{abstract}
We consider closed manifolds that admit a metric locally isometric to a product of symmetric planes. For such manifolds, we prove that the Euler characteristic is an obstruction to the existence of flat structures, confirming an old conjecture proved by Milnor in dimension $2$. In particular, the Chern conjecture follows in these cases. The proof goes via a new sharp Milnor--Wood inequality for Riemannian manifolds that are locally a product of hyperbolic planes. Furthermore, we analyze the possible flat vector bundles over such manifolds. Over closed Hilbert--Blumenthal modular varieties, we show that there are finitely many flat structures with nonzero Euler number and none of them corresponds to the tangent bundle.
Some of the main results were announced in C.R. Acad. Sci. Paris, Ser. I 346 (2008) 661-666.
\end{abstract}
\maketitle


\section{Introduction}

Let $M$ be a closed oriented smooth manifold. A principal $\mathrm{GL}^+(m,\mathbb{R})$-bundle $\xi$ over $M$ is called {\it flat} if it admits a flat structure, i.e.\ a connection on $\xi$ with zero curvature. Equivalently, a $\GL^+(m,\BR)$-bundle is flat if it is induced by a representation of the fundamental group $\pi_1(M)$. 
The manifold $M$ is said to admit a flat structure if the $\mathrm{GL}^+(m,\mathbb{R})$-bundle associated to its tangent bundle is flat, where $m$ now is the dimension of $M$. 

It is an old conjecture (see \cite{Mi58,KaTo68,HiTh75}) that the Euler characteristic is an obstruction to the existence of flat structure: 

\begin{conj}\label{conj}
A closed aspherical manifold with nonzero Euler characteristic does not admit a flat structure.
\end{conj}

Originally the conjecture was stated for general closed manifolds and not just for aspherical ones, but a counterexample was found by Smillie \cite{Sm77}. 

Conjecture \ref{conj} was proved for surfaces by Milnor in his celebrated paper \cite{Mi58}. For dimension greater than 2 however, very little progress has been made. Hirsch and Thurston \cite{HiTh75} confirmed it for manifolds whose fundamental group is a free product of virtually solvable groups\footnote{In their 1975 paper, Hirsch and Thurston allowed also free factors of polynomial growth. However, by Gromov's \cite{Gr81} theorem, such groups are virtually nilpotent and in particular virtually solvable.}.
One of our main results is the following:

\begin{thm}\label{thm:flat}
Let $M$ be a closed oriented manifold admitting a Riemannian structure with respect to which the universal cover is a direct product of $2$-dimensional symmetric spaces. If $\chi(M)\ne 0$ then $M$ admits no flat structure.
\end{thm}

An important example of flat structures on manifolds are affine structures, and the famous Chern conjecture that closed affine manifolds have zero Euler characteristic is a particular case of Conjecture \ref{conj}. As far as we know, Chern's conjecture might also hold for nonaspherical manifolds. Note that since the Euler characteristic vanishes in odd dimensions, these conjectures concern even dimensional manifolds only.  

It follows in particular that the Chern conjecture holds for manifolds as in Theorem \ref{thm:flat}. Note that in the case of products (even of hyperbolic surfaces) the nonexistence of affine structures on the product cannot be directly deduced from the nonexistence of affine structure on the factors (cf. Example \ref{Example Ghys} below).
In contrary to Conjecture \ref{conj} for which there were hardly any previous results in dimension $>2$ apart from \cite{HiTh75}, Chern's conjecture was confirmed in a few additional special cases: Goldman and Hirsch \cite{GoHi84} proved that higher rank {\it irreducible} locally symmetric manifolds can never admit an affine structure. Moreover, Kostant and Sullivan \cite{KoSu75} proved that a closed manifold with nonzero Euler characteristic cannot admit a {\it complete} affine structure. However, proving the nonexistence of a {\it non}-complete affine structure is usually much harder. 


In general, manifolds as in Theorem \ref{thm:flat} with vanishing Euler characteristic may admit a flat and even an affine structure. For instance $\gS_2\times S^2\times T^2$, where $\gS_2$ is a surface of genus $2$, $S^2$ the $2$-sphere and $T^2$ is the $2$-dimensional torus, admits an affine structure. More generally, suppose $M$ is a compact Riemannian manifold whose universal cover is a product of $2$-dimensional symmetric spaces, i.e. 
$\ti M\cong \mathcal{H}^{n_1}\times (S^2)^{n_2}\times \BR^{2n_3}$, where $\mathcal{H}$ is the hyperbolic plane. If $n_3\ne 0$ we have $\chi(M)=0$ while if $n_2\ne 0$ we have that some finite cover of $M$ is an $S^2$-bundle over a smaller dimensional closed manifold, in which case the Euler class of any flat oriented $\BR^{\dim(M)}$-bundle over $M$ vanishes. Thus the general case reduces to the case $n_2=n_3=0$. 
In most of the paper we will restrict to this case. Theorem \ref{thm:flat} reduces to Corollary \ref{cor} which is an immediate consequence of the sharp Milnor-Wood inequality, Theorem \ref{Theorem: Milnor-Wood for products}.

Although we consider manifolds of a very special kind, we believe that some of the ideas are extendable for a broader setup and that more general cases of Conjecture \ref{conj} will be settled in the future. One may hope to prove it for more general locally symmetric spaces or Hadamard manifolds.


\subsection*{The Euler class}
When studying flat oriented vector bundles using characteristic classes one should focus on the Euler class - the only class that carries a nontrivial data (cf. Remark \ref{rem:ECnot0}).

Let $\xi$ be a principal $\mathrm{GL}^+(m,\mathbb{R})$-bundle.
The (real) {\it Euler class} $\varepsilon_m(\xi)\in H^m(M,\mathbb{R})$ of $\xi$ is the image under the inclusion of coefficients $\mathbb{Z}\hookrightarrow\mathbb{R}$ 
of the Poincar\'e dual of the zero locus of a generic section in the associated oriented $\BR^m$ bundle (the first obstruction to the existence of a nowhere vanishing section). The {\it Euler number} of $\xi$ is the natural pairing of the Euler class with the (real) fundamental class $[M]\in H_m(M,\mathbb{R})$ of $M$:
$$\chi(\xi)=\langle \varepsilon_m(\xi),[M]\rangle.$$

Lusztig and Sullivan \cite{Su76} observed that there are only finitely many
isomorphy classes of $\mathrm{GL}^+(m,\mathbb{R})$-bundles admitting a flat structure, and  hence a bound depending only on $M$ for the possible Euler numbers of flat $\mathrm{GL}^+(m,\mathbb{R})$-bundles over $M$. Indeed, two representations in the same connected component of the space of representations $\text{Rep}(\pi_1(M),\GL^+(m,\BR))$ induce isomorphic bundles, and the space of representations is a real algebraic variety, hence it has finitely many connected components.\footnote{Note that for a given bundle there could still be uncountably many flat structures -- corresponding to nonconjugate representations lying in the same component.}

Milnor \cite{Mi58}  proved that a $\GL^+(2,\mathbb{R})$-bundle $\xi$ over a surface $\Sigma_g$ of
genus $g\geq 1$ admits a flat structure if and only if its Euler
number $\chi(\xi)$ satisfies the inequality
$|\chi(\xi)|
\leq g-1$.
In particular this shows that contrary to the Chern and Pontrjagin classes, the Euler class is a nontrivial characteristic class for flat bundles. Indeed, every integer can be realized as an Euler number of some $\GL^+(2,\BR)$-bundle over a surface. 
Milnor's inequality was later generalized to circle bundles by Wood
\cite{Wo75}.

In dimension greater than $2$, there were up to now few examples where explicit bounds were given for $\chi(\xi)$, or more generally for any primary characteristic number. In an unpublished work, Smillie established explicit bounds for the Euler number of flat $\GL^+(2n,\mathbb{R})$-bundles $\xi$ over hyperbolic manifolds $M$ of even dimension $2n$:
\begin{eqnarray}
 |\chi(\xi)| \leq \frac{\pi^n}{2^n\cdot 1\cdot 3\cdot 5 \cdot \ldots \cdot (2n-1) v_{2n}} | \chi(M) |.\label{eqn: Smillie}
\end{eqnarray}
However, it is not known whether nontrivial flat bundles over such manifolds exist at all when $n>1$. Moreover, observe that the constant is strictly greater than $1$ when $n>1$, so that this inequality does not imply the nonexistence of flat or affine structures on $M$. In a more geometric direction, Besson, Courtois and Gallot  \cite{BCG} proved sharp Milnor--Wood type inequalities for the pullback of the volume form under
representations $\Gamma \rightarrow \mathrm{Isom}_+(X)$ of cocompact lattices $\Gamma$ in $\mathrm{Isom}_+(X)$, where $X$ is a product of symmetric spaces of strictly negative curvature. Note however that in dimension $>2$, the volume form is in general not a primary characteristic class. One purpose of the current paper is to prove sharp generalizations of Milnor's inequality for a family of higher dimensional manifolds and use these inequalities to analyze the possible flat bundles and confirm in particular Conjecture \ref{conj}.


\subsection*{Statement of the main results}

Denote by $\mathcal{H}$ the real hyperbolic plane. We consider closed manifolds $M$ admitting a complete Riemannian structure locally isometric to $\mathcal{H}^n$ for some $n\geq1$. 
It is well known that the universal cover of $M$ is isometric to the symmetric space $\mathcal{H}^n$ and the fundamental group $\pi_1(M)$ acts on $\mathcal{H}^n$ by deck transformations. This produces an embedding of $\pi_1(M)$ as a torsion free cocompact lattice in $\text{Isom}(\mathcal{H}^n)$.
We call such $M$ an $\mathcal{H}^n$-manifold.
We prove Milnor--Wood type inequalities for $\mathcal{H}^n$-manifolds:

\begin{thm}
\label{Theorem: Milnor-Wood for products}
Let $M$ be a closed $\mathcal{H}^n$%
-manifold and $\xi$ a
$\mathrm{GL}^+(2n,\mathbb{R})$-bundle over $M$. If $\xi$ admits a
flat structure, then
\[
| \chi(\xi) |=| \left\langle \varepsilon(\xi),[M]\right\rangle |
\leq\frac{1}{\left( -2\right)  ^{n}}\chi(M).
\footnote{Note that the right hand side of the inequality is always strictly
positive, since the Euler characteristic $\chi(M)$ is nonzero, and
its sign is $(-1)^n$. This can for example be seen from Hirzebruch's
Proportionality Principle \cite{Hi58} recalled below.}
\]
\end{thm}

The case $n=1$ is Milnor's celebrated inequality. 
For $n=2$, a weaker upper bound is obtained in \cite{Bu07} by combining
the explicit computation of the simplicial volume
for such manifolds and Ivanov--Turaev's upper bound
\cite{IvTu82} for the Euler class of flat oriented vector bundles.

\begin{cor}\label{cor}
A closed $\mathcal{H}^n$-manifold does not admit a flat structure.
\end{cor}

This result is new for $n\geq3$. For $n=1$ and $n=2$ it follows from the corresponding inequalities of \cite{Mi58} and \cite{Bu07} respectively.  

\begin{proof}
[Proof of Corollary \ref{cor}]
By definition, $M$ admits a flat structure if and only if its tangent bundle $TM$ admits a flat structure. Thus by Theorem \ref{Theorem: Milnor-Wood for products},%
\[
|\chi(M)|=|\chi(TM)|=| \left\langle
\varepsilon(\xi),[M]\right\rangle | \leq\frac {1}{2^{n}}|\chi(M)|,
\]
which is impossible since $\chi(M)\neq 0$.
\end{proof}

Note that Theorem \ref{thm:flat} follows from Corollary \ref{cor}. For general Riemannian manifolds as in Theorem \ref{thm:flat} whose universal cover has a $2$-sphere or a Euclidean plane as a factor the Euler number of a flat oriented vector bundle always vanishes, and the corresponding Milnor--Wood inequality is consequently trivial.

The $\mathcal{H}^n$-manifolds are of particular interest among all locally symmetric manifolds. Indeed, while one can deduce from superrigidity theorems that some rigid locally symmetric manifolds $M$ admit no nontrivial flat bundle of dimension $\dim (M)$, $\mathcal{H}^n$-manifolds do admit (in many cases a unique) flat bundle with nonzero Euler number (see Theorems \ref{thm: gnrl Euler number} and \ref{thm:rigidity} and Corollaries \ref{cor: Euler number for rigid} and \ref{cor:unique}). By Corollary \ref{cor} these bundles cannot be isomorphic to the tangent bundle $TM$. 

\medskip

Note that $\text{Isom}(\mathcal{H}^n)\cong S_n\ltimes\PGL(2,\BR)^n$ where $S_n$ denotes the full $n$-th symmetric group. The orientation-preserving isometries form a subgroup of index $2$ denoted $\text{Isom}^+(\mathcal{H}^n)$. Set $G_n=S_n\ltimes\prod_{i=1}^n\GL^1(2,\BR)$, where $\GL^1(2,\BR)$ denotes the group of $2\times 2$ real matrices with determinant $\pm 1$. Then $G_n$ is an order $2^n$ sheeted cover of $\text{Isom}(\mathcal{H}^n)$. Let $G_n^+$ be the preimage of $\text{Isom}^+(\mathcal{H}^n)$ in $G_n$. Note that $G_n$ admits a natural faithful $2n$-dimensional real representation, where the image of $G_n^+$ is the intersection of $G_n$ with $\SL(2n,\BR)$. We denote by $s:\text{Isom}^+(\mathcal{H}^n)\to S_n$ the canonical (surjective) homomorphism, and by $f:G_n^+\to\text{Isom}^+(\mathcal{H}^n)$ the covering map.

\begin{defn}
A discrete subgroup $\gC\le\text{Isom}^+(\mathcal{H}^n)$ will be called {\it cofaithful} if it admits a faithful lift to $G_n^+$, i.e. if there exists a subgroup $\ti\gC\le G_n^+$ which intersects $\ker(f)$  trivially and such that $\gC=f (\ti\gC)$. 
In that case, the isomorphism $\gC\to\ti\gC$ inverse to the restriction of $f$ to $\ti\gC$ will be called a {\it cofaithful lift} (or a cofaithful map). An $\mathcal{H}^n$-manifold $M$ will be called {\it cofaithful} if its fundamental group is cofaithful.
\end{defn}

\begin{rem}
Milnor proved \cite{Mi58} that every compact hyperbolic surface is cofaithful. Therefore a direct product of surfaces is also cofaithful. However, we do not know if for $n>1$ every compact $\mathcal{H}^n$-manifold is cofaithful.

If $\ti\gC\le G_n^+$ has no elements of order two, then $\gC=f(\ti\gC)$ is cofaithful. Since every finitely generated subgroup of $G_n^+$ admits a subgroup of finite index with this property (for instance by \cite[6.11]{Raghunathan}), it follows that any complete $\mathcal{H}^n$-manifold admits a cofaithful finite cover. 
\end{rem}

The next result gives the precise value of the Euler number of flat vector bundles induced by a cofaithful map.

\begin{thm}\label{Theorem: sharp}
The flat 
$\mathrm{GL}^+(2n,\mathbb{R})$-bundle $\xi$ 
on a closed cofaithful orientable manifold $M=\gC\backslash \mathcal{H}^n$
induced by a cofaithful map $\gC\to G_n^+$ composed with the natural embedding of $G_n^+$ in $\GL^+_{2n}(\BR)$ satisfies
\[
| \chi(\xi) |=\frac{1}{2^{n}}|\chi(M)|.
\]
\end{thm}

Thus, up to a finite cover, the inequality of Theorem \ref{Theorem: Milnor-Wood for products} is sharp for every closed $\mathcal{H}^n$-manifold. 
However, not every integer in the interval
$[\frac{-|\chi(M)|}{2^n},\frac{|\chi(M)|}{2^n}]$ is in general an
Euler number of a flat bundle, and Theorem \ref{Theorem: Milnor-Wood
for products} can be refined as in Theorem \ref{thm: gnrl Euler
number} below. 

We will say that an $\mathcal{H}^n$-manifold is \textit{rigid} if it
has no finite cover which decomposes as a product manifold with a
2-dimensional factor. This terminology is motivated by the (local,
Mostow and Margulis) rigidity theorems which apply for such
manifolds. For instance, Hilbert--Blumenthal modular varieties are rigid. By a closed Hilbert--Blumenthal\footnote{Note that the varieties considered in the original work of Blumenthal and Hilbert were isotropic, i.e. noncompact while here the unisotropic (or closed) ones are considered. These are sometimes referred as quaternionic Shimura varieties or Picard modular manifolds.} manifold we mean a manifold of the form $\gC\backslash\mathcal{H}^n$ where $\gC$ is an anisotropic Hilbert--Blumenthal group, i.e. a subgroup of finite index of $\BG(\OO_\BK)$ where $\BK$ is a totally real number field of degree $m\ge n>1$, $\OO_\BK$ its ring of integers and $\BG$ an anisotropic $\BK$ almost simple algebraic group such that $\prod_{\gc\in\text{Gal}(\BK/\BQ)}\BG^\gc(\BR)\cong\SL_2(\BR)^n\times K$ for some compact group $K$. The induced embedding of $\gC$ in $\SL_2(\BR)^n$ is discrete cocompact and irreducible.  
By Margulis' arithmeticity theorems, every rigid $\mathcal{H}^n$-manifold admits a finite cover which is a product of Hilbert--Blumenthal modular manifolds.

Recall (see \cite[Theorem 5.22]{Raghunathan}) that every closed $\mathcal{H}^n$-manifold $M$ admits a finite cover $N$ of the form 
$$N=\Sigma_{g_1}\times\ldots\times \Sigma_{g_k}\times N',
$$
where $N'$ is rigid, $k\ge 0$, and the $\Sigma_{g_i}$'s are surfaces of genus $g_i\ge 2$. The pullback of a flat $\GL^+(2n,\mathbb{R})$-bundle $\xi_M$ over $M$ is a flat $\GL^+(2n,\mathbb{R})$-bundle over $N$.

\begin{thm} \label{thm: gnrl Euler number}
Let $N$ be as above, and let $\xi_N$ be a flat $\GL^+(2n,\mathbb{R})$-bundle over $N$. Then
$$
 \chi(\xi_N)\in \bigg\{\pm \frac{\chi(N')}{2^{{\mathrm{dim}(N')}/{2}}}\prod_{i=1}^k \ell_i :  |\ell_i|\leq g_i-1\bigg\}\cup\{0\}.
$$
Moreover, if $N$ is cofaithful, all the integers are actually attained as Euler numbers of flat bundles.
\end{thm}

Rewriting Theorem \ref{thm: gnrl Euler number} in two special interesting cases, namely when $M$ is a product of surfaces, or when $M$ is rigid, we obtain:

\begin{cor}\label{cor: Euler number of surfaces}
Let $\xi$ be a $\mathrm{GL}^+(2n,\mathbb{R})$-bundle over the
product of closed oriented surfaces $M_{1},\ldots,M_{n}$ of respective
genus $g_i\geq 2$, for $i=1,\ldots,n$. If $\xi$ is flat, then
\[
\chi(\xi)\in\bigg\{
\prod_{i=1}^n
\ell_{i}: \text{ }\left\vert \ell_{i}\right\vert \leq g_i-1\bigg\} .
\]
\end{cor}

\begin{cor}\label{cor: Euler number for rigid} 
Let $M$ be a compact rigid $\mathcal{H}^n$-manifold, and let $\xi$ be a
$\mathrm{GL}^+(2n,\mathbb{R})$-bundle over $M$. If $\xi$ is flat
and $\chi(\xi)\neq 0$, then
$$|\chi(\xi)|=\bigg|\frac{1}{2^n}\chi(M)\bigg|.$$
\end{cor}


For rigid manifolds we show that the flat $\dim(M)$-dimensional vector bundles with positive Euler number are parameterized by the (finite) first
cohomology group in the $n$-dimensional vector space over the field $\BF_2$ of two elements. More precisely, let $M$ be a closed oriented $\mathcal{H}^n$-manifold with fundamental group $\gC\le\text{Isom}(\mathcal{H}^n)$. We denote by $H^1_\text{\rm Sym}(M,\BF_2^n)=H^1_\text{\rm Sym}(\gC,\BF_2^n)$ the first cohomology group of $\gC$ in $\BF_2^n$ with respect to the action of $\gC$ on $\BF_2^n$ that permutes the coordinates according to the map $s:\gS\to S_n$. 

\begin{thm}\label{thm:rigidity}
Let $M$ be a compact rigid $\mathcal{H}^n$-manifold. Then there exists a flat vector bundle over $M$ with nonzero Euler number if and only if $M$ is cofaithful. In that case, there is a one-to-one correspondence between the set of flat vector bundles over $M$ with positive Euler number and the elements of the finite group $H^1_\text{\rm Sym}(M,\BF_2^n)$. Moreover $M$ admits a canonical normal finite cover $N$ such that all these bundles induce the same vector bundle with the same flat structure over $N$.  
\end{thm}

When $M$ is a Hilbert--Blumenthal modular manifold or a direct product of such, the fundamental group $\gC$ is embedded in the direct product $\prod_{i=1}^n\PSL(2,\BR)$ and in particular the map $s:\gC\to S_n$ is trivial. In that case $H^1_\text{Sym}(M,\BF_2^n)\cong \text{Hom}(\gC,\BZ/(2))^n$. Hence we have:

\begin{cor}
Let $M$ be a closed cofaithful Hilbert--Blumenthal modular manifold. Then $M$ admits exactly $2^{nk+1}$ flat $2n$-dimensional vector bundles with nonzero Euler number, where $k=\dim H^1(M,\BZ/(2))$ and $n=\dim (M)/2$. The Euler number for each of these vector bundles is either $\pm\chi (M)/2^n$, where the sign is determined according to the choice of orientation.
\end{cor} 

Note that $\dim H^1(M,\BZ/(2))$ is always finite since $\pi_1(M)$ is finitely generated.
On the other hand, there are many examples of closed Hilbert--Blumenthal modular manifolds with nontrivial first $\BZ/(2)$-cohomology. In fact if $\gC\le\prod_{i=1}^n\SL_2(\BR)$ is an arithmetic lattice, one can map it to $\SL_2(\BF_q)$ for some finite field $\BF_q$, and take $\gC'$ to be the preimage of the $2$-Sylow subgroup. Then $H^1(\gC',\BZ/(2))\ne\langle 0\rangle$. More generally, Lubotzky \cite{Lubotzky} showed that every finitely generated linear group admits a finite index subgroup which is mapped onto $\BZ/(2)$. When $k=0$ we have:\footnote{This uniqueness result was stated incorrectly in \cite{BuGe08}, where the condition $k=0$ is missing.}

\begin{cor}\label{cor:unique} 
Let $M$ be a compact cofaithful Hilbert--Blumenthal modular $\mathcal{H}^n$-manifold with trivial first $\BF_2$-cohomology.
Then there exists a unique flat $\mathrm{GL}^+(2n,\mathbb{R})$-bundle over $M$ with positive Euler number.
\end{cor} 

In contrast to Milnor's characterization of flat bundles over
surfaces \cite{Mi58} which states that a $\GL^+(2,\bbr)$-bundle
$\xi$ over a surface of genus $g$ is flat {\it if and only if}
$|\chi(\xi)|\leq g-1$, the converse of Theorems \ref{Theorem:
Milnor-Wood for products} and \ref{thm: gnrl Euler number} do not hold in general;
there exist nonflat bundles whose Euler numbers are allowed as possible Euler
numbers of flat bundles. Examples are given
in Section \ref{section:Examples}.

\begin{rem} Our results also hold for unoriented manifolds. Indeed, the Euler number of an oriented vector bundle can also be defined over an unoriented manifold $M$, for example as one half of the Euler number of the pullback of the bundle to an oriented double cover of $M$. Similarly, the simplicial volume of an unoriented manifold is simply one half of the simplicial volume of an oriented double cover. All our proofs apply with minor changes to this wider setting. 
\end{rem}


\subsection*{Outline of the proof of Theorem \ref{Theorem: Milnor-Wood for products}}

In order to guide the reader through this paper and the proof of the generalized Milnor--Wood inequality,  to which most of the paper is dedicated, let us describe the skeleton of the proof of Theorem \ref{Theorem: Milnor-Wood for products}.

Let $M$ be a closed $\mathcal{H}^n$-manifold and let $\xi$ be an oriented $2n$-vector bundle over $M$. Assume that $\xi$ is flat, then it induces a representation $\rho:\pi_1(M)\to\GL^+(2n,\BR)$, and $\xi$ is isomorphic to 
$$
 \pi_1(M)\backslash (\mathcal{H}^n\times\BR^{2n}),
$$
where $\pi_1(M)$ acts on $\mathcal{H}^n$ by deck transformations and on $\BR^{2n}$ via $\rho$.

As explained in Section \ref{section: the Euler class}, the Euler class 
$\varepsilon_{m}(\xi)$ is the pullback under $\rho$ of the Euler class $\varepsilon_{m}\in H_{c}^{m}(\GL^{+}(m,\bbr))$ which is continuous and admits a representing cocycle of bounded $L^\infty$ norm.

By replacing $M$ with a finite cover, which amounts to multiplying both sizes of the inequality in Theorem \ref{Theorem: Milnor-Wood for products} by the degree of the cover, we may assume that the image $\gC$ of $\pi_1(M)$ in $\text{Isom}(\mathcal{H}^n)$ lies in the identity connected component $\text{Isom}(\mathcal{H}^n)^\circ\cong\prod_{i=1}^n\PSL(2,\BR)$. Denoting by $G_i$ the closure of the projection of $\gC$ to the $i$'th factor of $\prod_{i=1}^n\PSL(2,\BR)$, and using Margulis superrigidity theorem we can extend (up to moding out some finite central subgroup of the image) $\rho$ to a representation $\ti\rho:\prod_{i=1}^nG_i\to \GL(2n,\BR)$. Hence we are led to analyze some properties of representations of product groups. A special feature of the case we consider is that the number of factors is exactly half the dimension of the representation. In a sense, this ``lucky situation" allows interesting examples of flat vector bundles to exist, but leaves no room for uncontrollable mysterious ones. 

For any representation $\ti\rho:\prod_{i=1}^nG_i\to \GL(2n,\BR)$ we associate a canonical integer $t(\ti\rho)$ (see Section \ref{Section: Representations of products}) which, vaguely speaking, measures how many tensor sub-representations are encoded in $\ti\rho$. In particular, $t(\ti\rho)\ne 0$ if and only if the restriction of $\ti\rho$ to the product of two of the associated factors $G_i\times G_j$ and to some subspace of $\BR^n$ is isomorphic to the tensor product of some irreducible nontrivial representations of $G_i$ and $G_j$. 
Furthermore, if $t(\ti\rho)\ne 0$, switching the order of $i$ and $j$ produces a symmetry of the corresponding vector bundle which, in case the corresponding subrepresentations are both $2$ dimensional, reflects the orientation without changing the sign of the generator of the top dimensional cohomology. In particular this forces the Euler class to vanish
(see Lemmas \ref{lemma: tensor product} and \ref{lemma: tensor product factor}\label{lem:trho=0}).

We then prove the main Proposition \ref{prop:prod} which states, more or less, that there are $3$ (not necessarily exclusive) possibility for the image of $\ti\rho$. More precisely, up to dividing $\text{Image}(\ti\rho)$ by a normal amenable subgroup, which has no effect on norms of bounded cohomology classes (see Lemma \ref{lemma: can mod out amenable}), either
\begin{itemize}
\item $t(\ti\rho)\ne 0$, or
\item $\text{Image}(\ti\rho)$ acts faithfully by isometries on some symmetric space of dimension $< 2n$, or
\item $\text{Image}(\ti\rho)$ is conjugated to a subgroup of the block diagonal embedded $\prod_{i=1}^n\SL(2,\BR)\le\GL^+(2n,\BR)$.
\end{itemize}
In the first two cases $\rho^*(\varepsilon_{m})=0$, and we are left to study the third case. We deduce that the $L^\infty$-seminorm of the Euler class satisfies the inequality (see Theorem \ref{Thm: norm of the Euler class})
\begin{eqnarray}\label{Equ: Euler class inequality}
 \left\Vert \rho^{\ast}(\varepsilon_{2n})\right\Vert _{\infty}\leq
 \|\underbrace{\varepsilon_{2}\cup\ldots\cup\varepsilon_{2}}_{n\text{ times}}
 \|_\infty.
\end{eqnarray}

In parallel, using Hirzebruch's \cite{Hi58} and Gromov--Thurston's respective proportionality principles, we prove (see Proposition \ref{prop: Proportionality Principle Simpl Vol and Euler}) that the simplicial volume satisfies  
\begin{eqnarray}\label{Equ: prop princ for splcl vol in Intro}
 \left\Vert M\right\Vert =\frac{\chi(M)}{(-2)^{n}\left\Vert \varepsilon_{2}%
 \cup\ldots\cup\varepsilon_{2}\right\Vert _{\infty}}.
 \end{eqnarray}

Combining the last two equations we succeed to avoid the very hard calculations of the simplicial volume and the norm of the Euler class, and deduce that the Euler number satisfies
$$
 |\chi(\xi)|=| \left\langle \varepsilon(\xi),[M]\right\rangle| 
 \le\|\varepsilon(\xi)\|_\infty\cdot\| M \|
 \leq\frac{1}{\left( -2\right)  ^{n}}\chi(M).
$$

\medskip

In several places of the paper we will treat the rigid and the completely reducible case (when $M$ is a product of surfaces) separately, in order to obtain a better insight in each case.
In the rigid case, where $M$ is commensurable to a product of Hilbert--Blumenthal modular varieties, we show that $|\chi(\xi)|=0$ unless $M$ is cofaithful and $\rho$ is a twist of a cofaithful representation via a cocycle belonging to the finite cohomology group $H^1_\text{\rm Sym}(M,\BF_2^n)$ in which case $|\chi(\xi)|=\frac{1}{\left( -2\right)  ^{n}}\chi(M)$. 

\subsection*{Structure of the paper:} In Section \ref{section: Background} we discuss Milnor--Wood inequalities and their approach through bounded cohomology. After the historical introduction in Subsection \ref{subsection: MW inequalities via bdd coho}, we recall continuous bounded cohomology in Subsection \ref{Section: Simplicial volume and bounded cohomology} and specialize to the Euler class in Subsection \ref{section: the Euler class}. In Section \ref{section:Examples}, we give several examples of vector bundles, illustrating some issues regarding our main results. 
In Section \ref{section: Euler}, we prove vanishing results for the Euler class to be used in subsequent sections. In Section \ref{section: Proportionality Principles}, we establish the proportionality principle between the simplicial volume and the Euler characteristic, Equation (\ref{Equ: prop princ for splcl vol in Intro}). In Section \ref{Section: Representations of products} we study general product representations, prove a general result about representations in dimension twice the number of factors (Proposition \ref{prop:prod}) and define the integer $t(\rho)$ associated to a representation $\rho$. In Section \ref{Section: Representations}, we specialize to representations of lattices in $\PSL(2,\mathbb{R})^n$, applying the results of the previous section using Margulis superrigidity theorem. In Section \ref{Section:Estimating the norm of the Euler class}, combining the results of Sections \ref{section: Euler}, \ref{Section: Representations of products} and \ref{Section: Representations}, we prove Inequality (\ref{Equ: Euler class inequality}) for the Euler class of flat bundles. In Section \ref{section: Proof of the main theorems}, we prove the generalized Milnor--Wood inequality (Theorem \ref{Theorem: Milnor-Wood for products} and its refinement Theorem \ref{thm: gnrl Euler number}) and the fact that the inequality is indeed sharp (Theorem \ref{Theorem: sharp}). In Section \ref{section: pf of thm for rigid mflds} we discuss in more details the rigid case and prove Theorem \ref{thm:rigidity}, giving a complete characterization of flat bundles with a nonvanishing Euler number. 


\section{Background}\label{section: Background}
\subsection{Milnor--Wood inequalities via bounded cohomology}\label{subsection: MW inequalities via bdd coho} In his groundbreaking essay \cite{Gr82}, Gromov naturally puts
Milnor's inequality in the context of bounded cohomology. Indeed, canonical $L^1$ and $L^\infty$ norms can be defined on the spaces of
singular chains and cochains
. These in turn induce seminorms on the respective real
valued homologies and cohomologies. It then follows from the Hahn--Banach theorem that characteristic numbers can be expressed as the product of the corresponding seminorms. Let us be more precise:

Let $M$ be an $n$-dimensional closed oriented manifold. The $L^{1}$-norm
on the space $C_{\ast}(M)$ of real-valued chains on $M$, associated
to the canonical basis of singular simplices, 
\[
\left\Vert \sum_{i=1}^{r}a_{i}\sigma_{i}\right\Vert _{1}=\sum_{i=1}^{r}\left\vert a_{i}\right\vert ,
\]
for $\sum_{i=1}^{r}a_{i}\sigma_{i}$ in $C_{q}(M)$, induces a
seminorm, which we still denote by $\left\Vert -\right\Vert _{1}$,
on the real valued homology $H_{\ast}(M)$ of $M$. The seminorm of
a homology class is defined as the infimum of the norms of its representatives.
The \textit{simplicial volume} of $M$, denoted by $\left\Vert M\right\Vert $,
is the seminorm of the real valued fundamental class $\left[M\right]\in H_{n}(M)$
of $M$.

The dual $L^{\infty}$-norm (or \textit{Gromov norm}) on the space $C^{\ast}(M)$ 
of real valued cochains on $M$ is given, for every
cochain $c$ in $C^{q}(M)$, by\begin{eqnarray*}
\|c\|_{\infty} & = & \sup\{|c(z)|:z\in C_{q}(M)~\text{~with~}~\|z\|_{1}=1\}\\
 & = & \sup\{|c(\gs)|:\gs:\gD^{q}\to M~\text{~continuous}\}.\end{eqnarray*}
The subspace of bounded cochains $C_{b}^{\ast}(M)\subset C^{\ast}(M)$
consists of those cochains for which the Gromov norm is finite. Since
the boundary operators on $C_{\ast}(M)$ are bounded with respect
to the $L^{1}$-norm, the coboundary operators of $C^{\ast}(M)$
restrict to $C_{b}^{\ast}(M)$. The \textit{bounded cohomology} $H_{b}^{\ast}(M)$
of $M$ is, by definition, the cohomology of the cocomplex $C_{b}^{\ast}(M)$.
The inclusion of cocomplexes $C_{b}^{\ast}(M)\subset C^{\ast}(M)$
induces a \textit{comparison map} $c:H_{b}^{\ast}(M)\rightarrow H^{\ast}(M)$.
The Gromov norm on the space of (bounded) cochains induces a seminorm
on $H_{b}^{\ast}(M)$ and on $H^{\ast}(M)$, which we still denote
by $\left\Vert -\right\Vert _{\infty}$ (allowing the value $+\infty$
on $H^{\ast}(M)$). For $\alpha$ in $H^{q}(M)$, we have 
\[
 \left\Vert \alpha\right\Vert _{\infty}=\inf\{\|\alpha_{b}\|_{\infty}:\text{ }\alpha_{b}\in H_{b}^{q}(M),c(\alpha_{b})=\alpha\},
\]
where, over the empty set, the infimum is considered to be infinity.

It now follows from the Hahn--Banach theorem (see \cite[Corollary, page 7]{Gr82} or \cite[Proposition F.2.2]{BePe92}) that
\begin{equation}\label{Equ: b(M)=product of norms}
|\left\langle \beta,[M]\right\rangle|=\left\Vert
\beta\right\Vert_\infty\cdot\left\Vert M\right\Vert,~\forall\beta\in
H^n(M)~\mathrm{with}~\left\Vert \beta \right\Vert_\infty<\infty,
\end{equation}
where $\left\Vert M\right\Vert$ denotes the $L^1$ seminorm of the
fundamental class of $M$, the {\it simplicial volume} of
$M$.
Thus, if $\beta$ is a characteristic class, a bound on the
characteristic number $|\left\langle \beta ,[M]\right\rangle| $ can
be obtained by bounding both $\left\Vert \beta\right\Vert _{\infty}$
and $\left\Vert M\right\Vert $. Unfortunately, estimating each of
these terms is usually very difficult.

Nonzero exact simplicial volume computations are rare. For oriented
surfaces $\Sigma_{g}$ of genus $g\geq1$, it is not difficult to show
that $\left\Vert \Sigma_{g}\right\Vert =2\left\vert
\chi(\Sigma_{g})\right\vert =4(g-1)$. In particular, if $g\ge 2$ and
$\Sigma_{g}$ is endowed with a hyperbolic structure, then
$\left\Vert \Sigma_{g}\right\Vert =\
\mathrm{Vol}(\Sigma_{g})/\pi$. More generally, if $M$ is an
$m$-dimensional closed hyperbolic manifold, then $\left\Vert
M\right\Vert = \mathrm{Vol}(M)/v_m$ \cite{Gr82,Th78}, where
$v_{m}$ denotes the supremum of the volumes of geodesic simplices in 
$m$-dimensional hyperbolic space, and is known explicitly in low
dimensions only. The only further computation of a nonzero
simplicial volume is given in \cite{Bu08} for manifolds locally
isomorphic to the product of two copies of the hyperbolic plane. In
this case, one has $\left\Vert M\right\Vert =6\cdot\chi
(M)=3/(2\pi^{2})\cdot \mathrm{Vol}(M)$.

Gromov proved \cite{Gr82} (see also \cite{Bu07}) that characteristic classes of flat
$G$-bundles have finite $L^\infty$ seminorm when $G$ is a real
algebraic subgroup of $\mathrm{GL}(m,\mathbb{R})$, but actual upper
bounds for their norms are only known in special cases. For the
Euler class $\varepsilon_m$, Gromov \cite{Gr82} obtained from
Sullivan--Smillie's \cite{Su76} corresponding simplicial results that $\left\Vert
\varepsilon_{m}(\xi)\right\Vert _{\infty}\leq1/2^{m}$, whenever
$\xi$ is a $\GL^+(m,\mathbb{R})$-bundle admitting a flat structure.
Independently, Ivanov and Turaev \cite{IvTu82} exhibited an explicit
bounded cocycle representing the Euler class of flat bundles,
producing the same bound. In degree $2$, sharp upper bounds for the
K\"ahler class were computed by Domic and Toledo \cite{DoTo} in
terms of the rank of the associated symmetric space.
Clerc and \O rsted \cite{ClOr} later generalized this to include all
Hermitian symmetric spaces.

In view of the (im)possible seminorm computations, sharp
generalizations of Milnor's inequality were essentially carried
out in degree $2$ only. Note however that in dimension $2$, bounded cohomology not only naturally leads to Milnor--Wood type inequalities, but can further be used to study rigidity properties of representations of surface groups (see \cite{BuIoWi08}, \cite{BuIo07}, \cite{BuIoLaWi05}). In higher dimensions, Smillie's Inequality (\ref{eqn: Smillie}) for the Euler number of flat bundles over hyperbolic manifolds of even dimension $m=2n$ are obtained by using the upper bound $\left\Vert \varepsilon_{2n}%
(TM)\right\Vert _{\infty}\leq 1/2^{2n}$ for the Euler class and the value $\left\Vert M\right\Vert=\mathrm{Vol}(M)/v_{2n}=((2\pi)^n/(1\cdot 3\cdot 5 \cdot \ldots \cdot (2n-1) v_{2n}))| \chi(M) |$ for the simplicial volume. 

\subsection{Bounded group cohomology}\label{Section: Simplicial volume and bounded cohomology}
Recall that the continuous cohomology
$H_{c}^{\ast}(G)$ of $G$ is the cohomology of the cocomplex $C_{c}^{\ast}(G)^{G}$
endowed with its natural homogeneous coboundary operator $\delta$:
The space of continuous cochains is given as 
\[
 C_{c}^{q}(G)=\{c:G^{q+1}\longrightarrow\mathbb{R}:c\text{ is continuous}\},
\]
and $C_{c}^{q}(G)^{G}$ is the subspace of left $G$-invariant cochains,
where $G$ acts diagonally on $G^{q+1}$. The coboundary operator $\delta:C_{c}^{q}(G)^{G}\rightarrow C_{c}^{q+1}(G)^G$
is defined as 
\[
 \delta f(g_{0},\ldots,g_{q})=\sum_{i=0}^{q}(-1)^{q}f(g_{0},\ldots,\hat{g_{i}},\ldots,g_{q}),
\]
for $f\in C_{c}^{q}(G)^{G}$ and $(g_{0},\ldots,g_{q})\in G^{q+1}$.
For $c$ in $C_{c}^{q}(G)$, let \[
\|c\|_{\infty}=\sup\{|c(g_{0},\ldots,g_{q})|:(g_{0},\ldots,g_{q})\in G^{q+1}\}.\]
 Set \[
C_{c,b}^{q}(G)=\{c\in C_{c}^{q}(G):\|c\|_{\infty}<+\infty\},\]
 and let $C_{c,b}^{\ast}(G)^{G}$ be the cocomplex of continuous bounded
$G$-invariant cochains. Clearly, the coboundary operator restricts
to $C_{c,b}^{\ast}(G)^{G}$, and the \textit{continuous bounded cohomology}
$H_{c,b}^{\ast}(G)$ of $G$ is defined as the cohomology of this
cocomplex.

The inclusion of cocomplexes $C_{c,b}^{\ast}(G)^{G}\subset C_{c}^{\ast}(G)^{G}$
induces a \textit{comparison map} $c:H_{c,b}^{\ast}(G)\rightarrow H_{c}^{\ast}(G)$.
As in the singular case, the sup norm induces seminorms both on $H_{c,b}^{\ast}(G)$
and on $H_{c}^{\ast}(G)$ (again allowing the value $+\infty$ in
the latter case) and we have, for any $\alpha$ in $H_{c}^{q}(G)$,
that\[
\left\Vert \alpha\right\Vert _{\infty}=\inf\{\|\alpha_{b}\|_{\infty}:\alpha_{b}\in H_{b}^{q}(M),c(\alpha_{b})=\alpha\}.\]

If $\Gamma$ is a discrete group, then the continuity condition is
void and we omit the term {}``continuous\textquotedblright\ and
the subscript {}``c\textquotedblright\ in the corresponding terminology
and notations. Note that the group cohomology $H^{\ast}(\Gamma)$
is then nothing but the Eilenberg-MacLane cohomology of $\Gamma$
\cite[Chapter 4, paragraph 5]{McL63}.

Let $f:H\rightarrow G$ be a continuous homomorphism between topological
groups. The induced maps $f^{*}:H_{c}^{*}(G)\rightarrow H_{c}^{*}(H)$
and $f^{*}:H_{c,b}^{*}(G)\rightarrow H_{c,b}^{*}(H)$ on cohomology
are defined at the cochain level by precomposing any cochain $c:G^{q+1}\rightarrow\bbr$
with $(f,\ldots,f):H^{q+1}\rightarrow G^{q+1}$. The sup norm is not
increased at the cochain level, hence for any $\alpha$ in $H_{c}^{q}(G)$
or in $H_{c,b}^{q}(G)$, we have $\|f^{*}(\alpha)\|_{\infty}\leq\|\alpha\|_{\infty}$.

Let $\Gamma$ now be the fundamental group of a closed manifold $M$.
As for standard singular cohomology, the natural map $M\rightarrow B\Gamma$,
where $B\Gamma$ denotes the classifying space of $\Gamma$-bundles,
induces a natural map in bounded cohomology \[
H_{b}^{*}(\Gamma)\cong H_{b}^{*}(B\Gamma)\rightarrow H_{b}^{*}(M).\]
In contrast to the standard case, Gromov \cite[Section 3.1]{Gr82}
(see also \cite{Iv85}) proved the remarkable theorem that this map
is an isometric isomorphism. For aspherical manifolds (which
we will exclusively be dealing with) this theorem is easy to prove.
Indeed, in that case it is classical that there is an isomorphism
$H^{*}(\Gamma)\cong H^{*}(M)$ between the standard cohomology groups.
Furthermore, it is easy to exhibit explicit cochain maps $C^{*}(\Gamma)^{\Gamma}\rightarrow C^{*}(M)$
and $C^{*}(M)\rightarrow C^{*}(\Gamma)^{\Gamma}$ realising the isometric isomorphisms
$H^{*}(\Gamma)\cong H^{*}(M)$ and $H_{b}^{*}(\Gamma)\cong H_{b}^{*}(M)$.


\subsection{The Euler class}
\label{section: the Euler class}

The classical theorem of Hopf stating that the index of a nondegenerate
vector field over a manifold $M$ is equal to the Euler characteristic
$\chi(M)$ can be reformulated by saying that the Euler number of
the tangent bundle is equal to the Euler characteristic of $M$ \cite[Corollary 11.12]{MiSt79}:
\[
\chi(TM)=\chi(M).\]
(In fact, this is the origin of the name {}``Euler class''.) By
Poincar\'e duality, $\chi(M)=0$ when the dimension of $M$ is odd.
More generally, the real Euler class of an odd dimensional vector
bundle always vanishes \cite[Property 9.4]{MiSt79}.

The universal Euler class $\varepsilon_{m}$ lives in the $m$-th
cohomology group of the classifying space $B\mathrm{GL}^{+}(m,\bbr)$.
Denoting by $\mathrm{GL}^{+}(m,\bbr){}^{\delta}$ the group $\mathrm{GL}^{+}(m,\bbr)$
endowed with the discrete topology, recall that the classifying space
$B\mathrm{GL}^+(m,\bbr)^{\delta}$ classifies flat $\mathrm{GL}^{+}(m,\bbr)$-bundles
and that its cohomology is by definition isomorphic to the Eilenberg--Maclane
cohomology of its fundamental group $\mathrm{GL}^{+}(m,\bbr)$ \cite[Chapter 4, paragraph 5]{McL63}.

We shall further see that it is natural to consider the Euler class
as a continuous cohomology class in $H_{c}^{m}(\mathrm{GL}^{+}(m,\bbr))$.
More generally, let $G$ be a reductive Lie group and $K$ a maximal
compact subgroup in $G$. We have a commutative diagram 
\[
\xymatrix{H^{*}(BG)\ar[r]\ar[d] & H^{*}(BG^{\delta})\ar[d]^{\cong}\\
H_{c}^{*}(G)\ar@{^{(}->}[r] & H^{*}(G).}
\]
The left vertical arrow can be described as follows: On the one hand,
the embedding $K\hookrightarrow G$ is a homotopy equivalence and
hence $H^{*}(BG)\cong H^{*}(BK)$. Since the Chern--Weil homomorphism
is an isomorphism for compact groups \cite{Bott73}, $H^{*}(BK)$
is isomorphic to the ring of $\text{Ad}(K)$-invariant polynomials
on the Lie algebra of $K$. On the other hand, by the Van Est isomorphism
\cite{VE55} (see also \cite[Corollary 7.2]{Gui80}), the continuous
cohomology $H_{c}^{*}(G)$ is isomorphic to the $G$-invariant differential
forms on $G/K$. The quotient map $G\rightarrow G/K$ is a $K$-bundle.
This bundle is endowed with the Maurer--Cartan connection. Chern--Weil
theory now naturally assigns to any $\text{Ad}(K)$-invariant polynomial
a differential form on $G/K$, which is $G$-invariant by the $G$-equivariance of the
Maurer--Cartan connection. For more details
and explicit formulas for primary characteristic classes of flat bundles
viewed as continuous cohomology classes, see \cite{Du76}.

To see that the lower horizontal arrow is injective one can argue
as follows: Since $G$ is reductive, by Borel's Theorem \cite{Borel} it admits
a cocompact lattice $\gC$, and the map $H_{c}^{*}(G)\rightarrow H^{*}(G)$
is factorized by the restriction $H_{c}^{*}(G)\hookrightarrow H^{*}(\Gamma)$,
which is injective, having as left inverse the transfer map (given
by integration over a fundamental domain for $\Gamma\setminus G$).

For $G=\GL^{+}(m,\bbr)$ and $m$ even, the image of the Euler class
in any of the above cohomology groups, which we still denote by $\varepsilon_{m}$,
is nontrivial. Since a $\GL^{+}(m,\bbr)$-bundle $\xi$ over $M$
admitting a flat structure is induced from a representation $\rho:\pi_{1}(M)\rightarrow\GL^{+}(m,\bbr)$
of the fundamental group of $M$, the Euler class $\varepsilon_{m}(\xi)$
is just the image under the natural homomorphism $H^{m}(\pi_{1}(M))\rightarrow H^{m}(M)$
of $\rho^{*}(\varepsilon_{m})$, where $\rho^{*}:H_{c}^{m}(G)\rightarrow H^{m}(\pi_{1}(M))$
is the homomorphism induced by $\rho$. The fact that the Euler class
$\varepsilon_{m}(\xi)$ is the pullback of the continuous (bounded)
cohomology class $\varepsilon_{m}\in H_{c}^{m}(\GL^{+}(m,\bbr))$
is fundamental in our approach. \medskip{}

\begin{rem}\label{rem:ECnot0} The Euler class is the only characteristic
class (for $\GL^{+}(m,\bbr)$) surviving the passage to the classifying
space of flat bundles. All other generators of $H^{*}(B\GL^{+}(m,\bbr))$
are Pontrjagin classes and vanish on flat bundles since they are in
the image of the Chern-Weil homomorphism \cite[Chapter XII]{KoNo69}. For a topological proof of this classical consequence of Chern-Weil theory in the differential
context see \cite[Section 6]{KaTo68}. 
\end{rem}

Recall \cite{IvTu82} that the Euler class can be represented by a
bounded cocycle, i.e. it lies in the image of the comparison map $H_{c,b}^{m}(\mathrm{GL}^{+}(m,\bbr))\rightarrow H_{c}^{m}(\mathrm{GL}^{+}(m,\bbr))$.

Define $\text{Bl}_{2,n}^{+}$ to be the subgroup of $\GL^{+}(2n,\bbr)$
consisting of the diagonal embedding of the product of $n$ copies
of $\GL^{+}(2,\bbr)$, i.e. the image of
the injective homomorphism \[
\begin{array}{rcl}
\rho{}_{\Delta}:\prod_{i=1}^{n}\GL^{+}(2,\bbr) & \longrightarrow & \text{Bl}_{2,n}^{+}<\GL^{+}(2n,\bbr)\\
(A_{1},\ldots,A_{n}) & \longmapsto & \left(\begin{array}{ccc}
A_{1}\\
 & \ddots\\
 &  & A_{n}\end{array}\right).\end{array}\]
It follows from the Whitney product formula for the Euler class that
the induced cohomology map 
$$
 \rho_\Delta^*:H_{c}^{2n}(\GL^{+}(2n,\bbr))\rightarrow H_{c}^{2n}(\text{Bl}_{2,n}^{+})\cong H_{c}^{2n}(\prod_{i=1}^{n}\GL^{+}(2,\bbr))
$$
sends $\varepsilon_{2n}$ to the cup product $\varepsilon_{2}\cup\ldots\cup\varepsilon_{2}$
of $n$ copies of $\varepsilon_{2}$.

Let us end this preliminary section with two basic facts regarding $\varepsilon_m$:

Conjugation by an element $g\in \GL(m,\mathbb{R})$ induces an isomorphism on the underlying (universal) unoriented vector bundles, which preserves orientation if $\mathrm{det}(g)>0$ and reverses orientation if $\mathrm{det}(g)<0$. This yields:

\begin{lem}\label{lemma: conjugation =00003D mult by sign of det} 
For $g\in\GL(m,\mathbb{R})$
let $\rho_{g}:\GL^{+}(m,\mathbb{R})\rightarrow\GL^{+}(m,\mathbb{R})$
be the corresponding conjugation $\rho_{g}(A)=gAg^{-1}$, $A\in\GL^{+}(m,\mathbb{R})$.
Then \[
\rho^*_{g}(\varepsilon_{m})=\left\{ \begin{array}{rl}
\varepsilon_{m}\ \  & \mathrm{if}\det(g)>0,\\
-\varepsilon_{m}\ \  & \mathrm{if}\det(g)<0.\end{array}\right.\]
\end{lem}



Since the Euler class $\varepsilon_m\in H_c^m(\GL^+(m,\BR))$ lies in the image of the composition  
$$
 H^m_c(\PSL(m,\mathbb{R})) \longrightarrow H^m_c(\SL(m,\mathbb{R}))\longrightarrow H^m_c(\GL^+(m,\mathbb{R}))
$$
induced by the natural projections, for $\rho:\gC\to\GL^+(m,\BR)$, the Euler class $\rho^*(\varepsilon_m)$ is determined by the projection of $\rho$ to $\PSL(m,\BR)$. Similarly, we have:

\begin{lem} \label{lem: reps with same proj to products of PSL(2,R)} 
Let $\rho_1,\rho_2:\Gamma \rightarrow \prod_{i=1}^n \GL^+(2,\mathbb{R})$ be two representations 
whose projections 
to $\prod_{i=1}^n\PSL(2,\mathbb{R})$ coincide. Then
$$ \rho_1^*(\varepsilon_2 \cup \ldots \cup \varepsilon_2)=\rho^*_2(\varepsilon_2\cup \ldots \cup \varepsilon_2).$$
\end{lem}

\medskip


\section{Some simple examples}\label{section:Examples}
The purpose of this section is to answer several questions which naturally arise and to give examples illustrating some issues regarding the main results. 

\subsection{Nonflat bundles with even Euler number}
\begin{clm}\label{clm:3.1}
Let $M$ be a closed oriented manifold of even dimension $2n$, and let $k\in 2\mathbb{Z}$ be an arbitrary even integer. Then there exists an oriented $\mathbb{R}^{2n}$-vector bundle $\xi$  over $M$ with $\chi(\xi)=k$. 
\end{clm}

Since any manifold admits a degree $1$ mapping onto the sphere of the same dimension, obtained by sending an open disk $U$ in $M$ diffeomorphically onto $S^{2n}\setminus \{x_0\}$ and its complement $M\setminus U$ to some fixed point $x_0\in S^{2n}$, Claim \ref{clm:3.1} follows from its validity for $S^{2n}$. 
The tangent bundle $TS^{2n}$ over $S^{2n}$ has Euler number $\chi(TS^{2n})=2$, and an oriented vector bundle over $S^{2n}$ with Euler number equal to $2d$ is obtained by taking the pullback of $TS^{2n}$ by a self-map of degree $d$. 

Moreover for $n=1,2,4$ the assertion holds for odd $k\in \mathbb{Z}$ as well: 
Letting $F$ be the complex numbers, the quaternions or the octonions, the projective space $\BP(F^2)$ (which is a sphere of dimension $2,4$ or $8$) carries a canonical $F$-bundle (the total space consists of pairs $([V],v)$ where $[V]\in \BP(F^2)$ and $v\in V$) with Euler number equal to $1$. For more details about the $4$-dimensional case, see \cite[Lemma 20.9]{MiSt79}.

Since only finitely many isomorphy classes of bundles of a fixed dimension over a given manifold $M$ can admit a flat structure, most of the bundles constructed above will not be flat. We will now show that in certain cases, and in particular when the base is an $\mathcal{H}^n$-manifold, a nonflat vector bundle can attain every (or every even) integer as an Euler number.

\begin{prop}\label{prop: not flat dim 4} 
Let $M$ be a closed oriented $4$-dimensional manifold, and let $k\in \mathbb{Z}$ be an arbitrary integer. Then there exists an oriented $\mathbb{R}^{4}$-vector bundle $\xi$ over $M$ not admitting a flat structure and such that $\chi(\xi)=k$.
\end{prop}

\begin{proof} Suppose first that $k\neq 0$. Let $\xi_1$ be the oriented $\mathbb{R}^{4}$-vector bundle over $S^4$ underlying the canonical $F$-bundle over the projective space $\BP(F^2)\simeq S^4$ considered earlier, for $F$ the quaternions. It is proven in \cite[Lemma 20.9]{MiSt79} that $\chi(\xi_1)=1$ and the first real Pontrjagin class of $\xi_1$ is nonzero. As above, one constructs the bundle $\xi$ over $M$ by taking the pullback of $\xi_1$ by a map $f:M\rightarrow S^4$ of degree $k$. It follows that $\chi(\xi)=k$ and the first real Pontrjagin class of $\xi$ is nonzero, which implies that the bundle cannot admit a flat structure.
The case $k=0$ is treated similarly based on the existence of bundles over $S^4$ with vanishing Euler number and nonvanishing  first real Pontrjagin class shown in \cite[Lemma 20.10]{MiSt79}.
\end{proof}

\begin{lem}\label{lemma: not flat greater dim} Let $k,n\in \mathbb{Z}$ be integers and suppose that $n\geq 2$. Then there exists an oriented $\mathbb{R}^{2n}$-vector bundle $\xi$ over the product $\prod_{i=1}^n S^2$ of $n$ copies of the sphere $S^2$ such that $\xi$ does not admit a flat structure and $\chi(\xi)=k$. 
\end{lem}

\begin{proof} Let $\eta_1$ be the oriented $\mathbb{R}^{2}$-vector bundle over $S^2$ underlying the canonical $\mathbb{C}$-bundle over the projective space $\BP(\mathbb{C}^2)\simeq S^2$. It is easy to check that $\chi(\eta_1)=1$. As explained in the proof of Proposition \ref{prop: not flat dim 4}, there exists an oriented $\mathbb{R}^{4}$-vector bundle $\xi_k$ over $S^2\times S^2$ such that $\chi(\xi_k)=k$ and the first real Pontrjagin class of $\xi_k$ does not vanish: $p_1(\xi_k)\neq 0 \in H^4(S^2\times S^2,\BR)$. 
For $\xi=\xi_k \times \Pi_{i=3}^n \eta_1 $, we have 
$$
 \chi(\xi)=\chi(\xi_k)=k.
$$
Furthermore, since, for dimension reasons, all the Pontrjagin classes of $ \Pi_{i=3}^n \eta_1$ vanish, 
$$
 p_1(\xi)=p_1(\xi_k \times  \Pi_{i=3}^n \eta_1)=p_1(\xi_k)\neq 0 \in H^4(\Pi_{i=1}^n S^2,\BR),
$$ 
so that $\xi$ does not admit a flat structure.
\end{proof}


\begin{prop} 
Let $M$ be a closed oriented $\mathcal{H}^n$-manifold with $n\geq 2$. There exists a finite cover $N$ of $M$ such that for every even integer $k\in 2\mathbb{Z}$ there exists an oriented $\mathbb{R}^{2n}$-vector bundle $\xi$ over $N$ such that $\xi$ does not admit a flat structure and $\chi(\xi)=k$. 
\end{prop}

This implies in particular that upon passing to a finite cover, the converses of Theorems \ref{Theorem: Milnor-Wood for products} and \ref{thm: gnrl Euler number}  are wrong in dimension greater or equal to $4$. That is, in contrast to the case $n=1$, the Euler number not only does not determine the isomorphy class of an oriented $\mathbb{R}^{2n}$-vector bundle, it does not even determine whether a bundle admits a flat structure or not.

\begin{proof} 
For $n=2$, this is a consequence of Proposition \ref{prop: not flat dim 4} even without assuming that $k$ is even and without passing to a finite cover. Let us thus suppose that $n\geq 3$.

It is shown in \cite{Ok01} that there exists a finite cover $N$ of any closed oriented $\mathcal{H}^n$-manifold which admits a tangential map of nonzero degree to the dual compact symmetric space. That is, there exists a map $f:N\rightarrow \prod_{i=1}^n S^2$ which furthermore induces Matsushima's map on the corresponding singular real cohomology groups, and consequently has degree $\chi(N)/2^n$. Thus for $k\in (\chi(N)/2^{n})\BZ$ the result follows immediately from Lemma \ref{lemma: not flat greater dim}. 

In order to treat general $k\in 2\BZ$, start with $\xi_0$, the pullback by $f:N\rightarrow \prod_{i=1}^n S^2$ of an oriented $\BR^{2n}$-vector bundle over $\prod_{i=1}^n S^2$ with $\chi(\xi_0)=0$ and nonvanishing first Pontrjagin class, established in (the proof of) Lemma \ref{lemma: not flat greater dim}. Then modify $\xi_0$ inside a closed disk $U$ of $N$ as follows, in order to obtain a bundle with any given even Euler number: Let $V\subset U$ be a smaller disk contained in the interior of $U$. Let $x_0$ be an interior point in $V$ and $y_0$ be a point in $S^{2n}$, and consider the wedge $N \vee S^{2n}$ obtained by identifying $x_0$ and $y_0$. Let $\phi:N\rightarrow N\vee S^{2n}$ be defined as the identity on $N\setminus U$, a diffeomorphism $U\setminus V \rightarrow U\setminus \{x_0\}$ restricting to the identity on $\partial U$ and extending to the constant map $\partial V \rightarrow y_0$, and a map $V\rightarrow S^{2n}$ wrapping $V$ around the sphere and sending $\partial V$ to $y_0$. Let $\xi_k$ be the pullback through $\phi$ of the bundle over $N\vee S^{2n}$ obtained by gluing the bundle $\xi_0$ to some bundle with Euler number $k$ over the sphere. Clearly, $\chi(\xi_k)=k$. 

To see that the bundle $\xi_k$ cannot admit a flat structure, note that there exists a $4$-dimensional CW-complex $X$ in $N$ which has nonzero first Pontrjagin class $\xi_0$. Since $X$ has codimension at least $2$ in $N$, the closed disk $U$ can be chosen disjoint from $X$, so that the restrictions of $\xi_0$ and $\xi_k$ to $X$ will agree, and the first Pontrjagin class of $\xi_k$ being nontrivial on $X$ will be nonzero on $N$ (and in fact equal to $p_1(\xi_0)$ as a cohomology class on $N$, not only on $X$). In particular, $p_1(\xi_k)\neq 0 \in H^4(N,\mathbb{R})$, and the bundle $\xi_k$ cannot admit a flat structure.  
\end{proof}

\begin{rem}
In the case of a product of surfaces, examples of nonflat bundles with any given Euler number are easily constructed using the fact that there exists a degree $1$ map from any surface onto $S^2$. 
\end{rem}

\subsection{Flat bundles with zero Euler number} 

Let $M$ be an $\mathcal{H}^n$-manifold. For simplicity we will assume that its fundamental group is contained in the identity connected component $\text{Isom}(\mathcal{H}^n)^\circ\cong\prod_{i=1}^n\PSL(2,\BR)$ and admits a cofaithful representation $\rho:\pi_1(M)\rightarrow \prod_{i=1}^n \SL(2,\mathbb{R})$. Let $I$ be any subset of $\{1,2,\ldots,n\}$ and denote by $\mathrm{proj}_I:\prod_{i=1}^n \SL(2,\mathbb{R})\rightarrow \Pi_{i\in I} \SL(2,\mathbb{R})$ the canonical projection. Denote by $\rho_I:\pi_1(M)\rightarrow \SL(2n,\mathbb{R})$ the composition of $\rho$ with $\mathrm{proj}_I$ followed by the natural inclusion $\prod_{i\in I} \SL(2,\mathbb{R})\hookrightarrow \prod_{i=1}^n \SL(2,\mathbb{R})$ and the diagonal embedding $\prod_{i=1}^n \SL(2,\mathbb{R}) \hookrightarrow \SL(2n,\mathbb{R})$. Let $\xi_I$ be the corresponding flat bundle over $M$. For $I=\{1,2,\ldots,n\}$, we have $\chi(\xi_I)=2^{-n}| \chi(M) |$, while for $I\neq \{1,2,\ldots,n\}$, we have $\chi(\xi_I)=0$ since the representation $\rho_I$ commutes through a Lie group whose associated symmetric space has dimension strictly smaller than $2n$. 

For $I\neq J$ the representations $\rho_I$ and $\rho_J$ cannot be conjugated since we can chose a sequence $\gc_m\in\pi_1(M)$ such that $\rho_I(\gc_m)\to 1$ while $\rho_J(\gc_m)\to\infty$. 

For noncofaithful manifolds, examples can be constructed by taking linear representations of some finite quotient or, if $\pi_1(M)$ lies in $\text{Isom}(\mathcal{H}^n)^\circ$, considering its embedding in $\prod_{i=1}^n \PSL(2,\mathbb{R})$ as a cocompact lattice composed with the projections on factors followed with some linear (oriented) representation of $\PSL(2,\mathbb{R})$.


\subsection{The variety of flat structures}
For a given oriented closed manifold $M$, the space of flat $m$ ($=\dim M$) oriented vector bundles over $M$ coincides with the character variety $\text{Hom}(\pi_1(M),\GL^+(m,\BR))/G$ since two flat bundles are isomorphic as flat bundles if and only if the corresponding representations of $\pi_1(M)$ are conjugated. On the other hand two representations $\rho_1,\rho_2:\pi_1(M)\to\GL^+(m,\BR)$ that lie in the same connected component of $\text{Hom}(\pi_1(M),\GL^+(m,\BR))$ induce abstractly isomorphic bundles. 

If $M$ is a rigid $\mathcal{H}^n$-manifold, it can be deduced from Margulis' superrigidity theorem that two representations that lie in the same connected component are conjugate, and moreover that the character variety is finite. However, when $M$ is a nonrigid closed $\mathcal{H}^n$-manifold, connected components of $\text{Hom}(\pi_1(M),\GL^+(2n,\BR))$ in general have dimension larger than the corresponding conjugacy classes, and one deduces that there are continuously many nonequivalent flat structures over some fixed given vector bundle.


\subsection{Affine product of nonaffine manifolds}

\begin{exam}[E. Ghys]\label{Example Ghys}
There exist two closed manifolds $M_1,M_2$ whose product $M_1\times M_2$ admits an affine structure while neither $M_1$ nor $M_2$ do.
\end{exam}

Let $M$ be a compact quaternionic hyperbolic surface, i.e. an $8$-dimensional manifold of the form $\gC\backslash X$ where $\gC$ is a torsion free cocompact lattice of $\text{Sp}(2,1)$ and $X$ is the associated symmetric space. Take
$$
 M_1=(M\times S^1)~~~~~~~\text{and}~~~~~~~M_2=S^3.
$$ 

To give an affine structure on $(M\times S^1)\times S^3$ one can argue as follows: Let $Q$ denote the quaternions. Let $V$ be the open ''light cone" in $Q^{2,1}$ (i.e. $Q^3$ equipped with the standard $(2,1)$ form). $\text{Sp}(n,1)$ acts on $Q^3$ preserving this form and $\gC\backslash V$ is homeomorphic to $M\times Q^\ast$, where $Q^\ast$ denotes the space of nonzero quaternions and is homeomorphic to $\BR^+\times S^3$. 
Denote by $t$ the affine homothety of multiplying by $2$. Then $(\gC\times\langle t\rangle)\backslash V$ is homeomorphic to $M\times S^1\times S^3$ and it inherits the affine structure from $V$, as the $\gC\times\langle t\rangle$ action is affine.

Now $S^3$ is compact simply connected and hence cannot admit an affine structure. Consider $M\times S^1$. Suppose it admits an affine structure, and let $\rho:\gC\times\BZ\to\GL^+_9(\BR)$ be the associated representation. Since $\text{Sp}(2,1)$ has no nontrivial $9$-dimensional representations, it follows from the Corlette and Gromov--Schoen superrigidity theorems that $\rho(\gC)$ is finite. Hence up to replacing $M$ by a finite cover, we may assume that $\rho(\gC)$ is trivial, and so up to replacing $S^1$ by a finite cover, we have that either $\rho$ has trivial image, which is impossible since
$M\times S^1$ has no flat Riemannian structure (for instance by Bieberbach's theorem), or that the image of $\rho$ is infinite cyclic. The last possibility cannot hold as well; in fact the closed affine manifolds with cyclic holonomy were classified by Smillie (in his unpublished thesis \cite{Sm77a}) who showed that they are all Hopf manifolds.


\section{Vanishing results for the Euler class}\label{section: Euler}

In this section we assemble several statements indicating the vanishing of the Euler class for certain representations that comes out in our proof.


\subsection{Vanishing of the Euler class for tensor representations}
Identifying $\mathbb{R}^{2}\otimes\mathbb{R}^{2}$ with $\mathbb{R}^{4}$
we obtain a representation of $\GL(2,\bbr)\times\GL(2,\bbr)$ which
we call the tensor representation and denote by $\rho_{\otimes}$
(we will also use this symbol to denote the restriction of this representation
to subgroups of $\GL(2,\bbr)\times\GL(2,\bbr)$). Note that with respect
to the basis $\left\{ e_{1}\otimes e_{1},e_{1}\otimes e_{2},e_{2}\otimes e_{1},e_{2}\otimes e_{2}\right\} $
of $\mathbb{R}^{2}\otimes\mathbb{R}^{2}$ induced by the tensor of
standard basis vectors $\left\{ e_{1},e_{2}\right\} $ of $\mathbb{R}^{2}$,
$\rho_{\otimes}$ takes the form\[
\rho_{\otimes}(A,B)=A\otimes B=\left(\begin{array}{cccc}
a_{11}b_{11} & a_{11}b_{12} & a_{12}b_{11} & a_{12}b_{12}\\
a_{11}b_{21} & a_{11}b_{22} & a_{12}b_{21} & a_{12}b_{22}\\
a_{21}b_{11} & a_{21}b_{12} & a_{22}b_{11} & a_{22}b_{12}\\
a_{21}b_{21} & a_{21}b_{22} & a_{22}b_{21} & a_{22}b_{22}\end{array}\right),\]
where \[
A=\left(\begin{array}{cc}
a_{11} & a_{12}\\
a_{21} & a_{22}\end{array}\right)\mbox{ and }B=\left(\begin{array}{cc}
b_{11} & b_{12}\\
b_{21} & b_{22}\end{array}\right).\]

\begin{lem}\label{lemma: tensor product} The pullback of the Euler class vanishes
under the tensor representation $\rho_{\otimes}:\SL(2,\mathbb{R})\times\SL(2,\mathbb{R})\rightarrow\SL(4,\mathbb{R})$:
\[
\rho_{\otimes}^{*}(\varepsilon_{4})=0.\]\end{lem}

\begin{proof}The proof relies on the simple fact that switching the two factors
changes the sign of the orientation of the tensor product, and hence of the Euler class, while it does not change
the sign of the generator of $H_{c}^{4}(\SL(2,\mathbb{R})\times\SL(2,\mathbb{R}))\cong\mathbb{R}.$
We give the detailed proof for the convenience of the reader.

Let $\tau:\SL(2,\mathbb{R})\times\SL(2,\mathbb{R})\rightarrow\SL(2,\mathbb{R})\times\SL(2,\mathbb{R})$
denote the homomorphism permuting the two factors: $\tau(A,B)=(B,A),$
for $(A,B)\in\SL(2,\mathbb{R})\times\SL(2,\mathbb{R})$. Let $E_{23}$
denote the odd permutation \[
E_{23}=\left(\begin{array}{cccc}
1 & 0 & 0 & 0\\
0 & 0 & 1 & 0\\
0 & 1 & 0 & 0\\
0 & 0 & 0 & 1\end{array}\right),
\]
and let $\rho_{E_{23}}:\SL(4,\mathbb{R})\rightarrow\SL(4,\mathbb{R})$
denote the corresponding conjugation by $E_{23}$, \[
\rho_{E_{23}}(A)=E_{23}AE_{23}^{-1},\]
 for $A\in\SL(4,\mathbb{R})$. We have a commutative diagram \[
\xymatrix{\SL(2,\mathbb{R})\times\SL(2,\mathbb{R})\ar^{\rho_{\otimes}}[d]\ar^{\tau}[r] & \SL(2,\mathbb{R})\times\SL(2,\mathbb{R})\ar^{\rho_{\otimes}}[d]\\
\SL(4,\mathbb{R})\ar^{\rho_{E_{23}}}[r] & \SL(4,\mathbb{R}),}
\]
so that in particular 
\[
\tau^{*}(\rho_{\otimes}^{*}(\varepsilon_{4}))=\rho_{\otimes}^{*}(\rho_{E_{23}}^{*}(\varepsilon_{4})).
\]
But since $E_{23}$ has negative determinant, by Lemma \ref{lemma: conjugation =00003D mult by sign of det},
$\rho_{E_{23}}^{*}(\varepsilon_{4})=-\varepsilon_{4}$. Thus
we get, on the one hand, that \begin{equation}
\tau^{*}(\rho_{\otimes}^{*}(\varepsilon_{4}))=-\rho_{\otimes}^{*}(\varepsilon_{4}).\label{eq: equ 0 for lemma on tensor}\end{equation}
 On the other hand, letting $E_{(13)(24)}$ denote the even permutation
\[
E_{(13)(24)}=\left(\begin{array}{cccc}
0 & 0 & 1 & 0\\
0 & 0 & 0 & 1\\
1 & 0 & 0 & 0\\
0 & 1 & 0 & 0\end{array}\right),\]
 and $\rho_{E_{(13)(24)}}$ the corresponding conjugation by $E_{(13)(24)}$,
we also have a commutative diagram 
\[
\xymatrix{\SL(2,\mathbb{R})\times\SL(2,\mathbb{R})\ar^{\rho_{\Delta}}[d]\ar^{\tau}[r] & \SL(2,\mathbb{R})\times\SL(2,\mathbb{R})\ar^{\rho_{\Delta}}[d]\\
\SL(4,\mathbb{R})\ar^{\rho_{E_{(13)(24)}}}[r] & \SL(4,\mathbb{R}),}
\]
 where $\rho_{\Delta}:\SL(2,\mathbb{R})\times\SL(2,\mathbb{R})\rightarrow\SL(4,\mathbb{R})$
is the diagonal embedding. By Lemma
\ref{lemma: conjugation =00003D mult by sign of det}, we have 
\begin{equation}
\tau^{*}(\rho_{\Delta}^{*}(\varepsilon_{4}))=\rho_\Delta^*(\rho^*_{E_{(13)(24)}}(\varepsilon_4))=\rho_{\Delta}^{*}(\varepsilon_{4}),\label{equ 1 for lemma on tensor}
\end{equation}
since $E_{(13)(24)}$ has positive determinant. The Whitney product formula for the Euler class implies $\rho_{\Delta}^{*}(\varepsilon_{4})=\varepsilon_{2}\cup\varepsilon_{2}$,
so we can rewrite Equation (\ref{equ 1 for lemma on tensor}) as 
\begin{equation}
\tau^{*}(\varepsilon_{2}\cup\varepsilon_{2})=\varepsilon_{2}\cup\varepsilon_{2}.\label{equ 2 for lemma on tensor}
\end{equation}
 Finally, since  $H^{4}(\SL(2,\mathbb{R})\times\SL(2,\mathbb{R}))$
is one dimensional, generated by the cup product of Euler classes
$\varepsilon_{2}\cup\varepsilon_{2}$, there must exist $\lambda\in\mathbb{R}$
with $\rho_{\otimes}^{*}(\varepsilon_{4})=\lambda\cdot\varepsilon_{2}\cup\varepsilon_{2}$.
It then follows from Equations (\ref{eq: equ 0 for lemma on tensor})
and (\ref{equ 2 for lemma on tensor}) that \[
\lambda\cdot\varepsilon_{2}\cup\varepsilon_{2}=\rho_{\otimes}^{*}(\varepsilon_{4})=-\tau^{*}(\rho_{\otimes}^{*}(\varepsilon_{4}))=-\tau^{*}(\lambda\cdot\varepsilon_{2}\cup\varepsilon_{2})=-\lambda\cdot\varepsilon_{2}\cup\varepsilon_{2}.\]
Hence $\lambda=0$, and consequently $\rho_\otimes^*(\varepsilon_4)=0$.
\end{proof}


\subsection{Triviality of top cohomology of one factor, implies vanishing of the Euler class} 

\begin{lem} \label{lemma: cont coho of factor = 0}
Let $M_1,M_2$ be closed aspherical oriented manifolds of respective
dimensions $m_1,m_2$ and fundamental groups $\Gamma_1,\Gamma_2$. Set
$m=m_1+m_2$ and let 
$$\rho:\Gamma_1 \times \Gamma_2 \longrightarrow \GL^+(m,\mathbb{R})$$
be a representation. Let $G_1$ be the Zariski closure of
$\rho(\Gamma_1)$. If $H^{m_1}_c(G_1)=0$ then
$\rho^*(\varepsilon_m)=0\in H^m(\Gamma_1\times \Gamma_2)$.
\end{lem}

\begin{proof} Let $G_2$ be the Zariski closure of $\rho(\Gamma_2)$. 
Denote by $\varphi:G_1\times G_2 \rightarrow \GL^+(m,\bbr)$ the
homomorphism induced by the inclusions $G_i<\GL^+(m,\bbr)$, for
$i=1,2$. Note that $\varphi$ is well defined since $G_1$ and $G_2$
are by construction commuting subgroups of $\GL^+(m,\bbr)$, but it
is not necessarily injective since $G_1$ and $G_2$ may have a
nontrivial central intersection. The representation
$\rho$ naturally induces a homomorphism
$$
 \overline{\rho}:\Gamma_1\times \Gamma_2\longrightarrow G_1\times G_2,
$$
defined on the factors $\Gamma_i$, for $i=1,2$, as the restriction of $\rho$ to $\Gamma_i$. The original representation $\rho$ is now nothing else but the
composition
\begin{equation*}
\xymatrix{ \Gamma_1\times\Gamma_2
\ar@/^1.5pc/[rr]^{\rho}\ar[r]^{\overline{\rho}} & G_1\times G_2
\ar[r]^{\varphi} & \GL^+(m,\bbr).}
\end{equation*}
The
induced map
$\rho^*:H^m_c(\GL^+(m,\mathbb{R}))\rightarrow H^m(\Gamma_1\times\Gamma_2)$ thus factors through $H^m_c(G_1\times G_2)$, and the lemma will follow from:\\

\noindent {\bf Claim.} The induced map
$\overline{\rho}^*:H^m_c(G_1\times G_2)\rightarrow
H^m(\Gamma_1\times\Gamma_2)$ is zero (in degree $m$).

\medskip

To prove the claim, note that the two latter cohomology groups
satisfy a K\"unneth product formula: For $H^m_c(G_1\times G_2)$, it
follows from the Van Est isomorphism \cite{VE55} and the validity of
the K\"unneth formula for relative Lie algebra cohomology
\cite[1.3]{BoWa00}. For $H^m(\Gamma_1\times\Gamma_2)$, it holds
since by asphericality, the cohomology groups of $\Gamma_1\times
\Gamma_2$, $\Gamma_1$ and $\Gamma_2$ are isomorphic to the
cohomology groups of $M_1\times M_2$, $M_1$ and $M_2$ respectively,
and the latter cohomology groups satisfy the K\"unneth formula (see
for example \cite[Theorem A.6]{MiSt79}). Furthermore, since all
isomorphisms are natural, the map $\rho^*$ now becomes
$$\overline{\rho}^*=\rho_1^* \otimes \rho_2^*: \bigoplus\limits_{p+q=m} H^p_c(G_1)\otimes H^q_c(G_2) \longrightarrow  \bigoplus\limits_{p+q=m} H^p(\Gamma_1)\otimes H^q( \Gamma_2),$$
where $\rho_i:\Gamma_i\rightarrow G_i$, for $i=1,2$, denotes the
restriction of $\rho$ to $\Gamma_i$. Since $m=m_1+m_2$, the only
nonzero summand in the latter direct sum is the one corresponding to
$p=m_1,q=m_2$. Thus, the restriction of $\overline{\rho}^*$
to all the summands other then $H^{m_1}_c(G_1)\otimes H^{m_2}_c(G_2)$ is trivial. Hence, the assumption
$H^{m_1}_c(G_1)=0$ gives $\overline{\rho}^*=0$ as
claimed.
\end{proof}

\subsection{Euler class and amenable factors}

One advantage of continuous bounded cohomology is its blindness
to amenable factors:


\begin{lem}
\label{lemma: rho(Gamma) amenable}Let $M_{1},M_{2}$ be closed
aspherical oriented manifolds of respective dimension
$\dim(M_{1}),\dim(M_{2})\geq1$ and fundamental groups
$\Gamma_{1},\Gamma_{2}$. Set $m=\dim(M_{1})+\dim(M_{2})$ and let
\[
\rho:\Gamma_{1}\times\Gamma_{2}\longrightarrow\text{GL}^+(m,\mathbb{R})%
\]
be a representation. If either $\rho(\Gamma_{1})$ or
$\rho(\Gamma_{2})$ is amenable, then
\[
\rho^{\ast}(\varepsilon_{m})=0\in H^{m}(\Gamma_{1}\times\Gamma_{2}).
\]

\end{lem}

\begin{proof}
Suppose that $\rho(\Gamma_{1})$ is amenable. Denote by $G_{i}$, for
$i=1,2$, the Zariski closure of $\rho(\Gamma_{i})$. By the Tits
alternative, $\rho(\gC_1)$ being amenable, is virtually solvable,
thus its Zariski closure $G_{1}$ is also virtually solvable hence
amenable. Therefore the projection $G_1\times G_{2}\rightarrow
G_{2}$ induces, in virtue of \cite[Lemma 7.5.10]{Mo01}, an isometric
isomorphism $H^*_{c,b}(G_2)\cong H^*_{c,b}(G_1\times G_2)$ on the
respective continuous bounded cohomology groups.

As in the proof of Lemma \ref{lemma: cont coho of factor = 0}, let
$\varphi:G_1\times G_2\rightarrow \GL^+(m,\bbr)$ be the homomorphism
induced by the inclusion of the factors. Again, observe that
$\rho=\varphi\circ \overline{\rho}$, where
$\overline{\rho}:\Gamma_1\times \Gamma_2\rightarrow G_1\times G_2$
is as in the proof of Lemma \ref{lemma: cont coho of factor = 0}. Consider the commutative diagram
\begin{equation*}
\xymatrix{ H^*_c(\text{GL}^+(m,\mathbb{R}))
\ar@/^2pc/[rrr]^{\rho^*}\ar[r]^{\varphi^*} & H^*_c(G_1\times G_2)
\ar[rr]^{\overline{\rho}^*}
&&H^*(\Gamma_1\times \Gamma_2)& \\
&&H^*_c(G_2)\ar[rr] \ar@{_{(}->}[lu] &&
H^*(\Gamma_2)\ar@{_{(}->}[lu] \\
H^*_{c,b}(\text{GL}^+(m,\mathbb{R}))\ar[uu]_c \ar[r]^{\varphi^*} &
H^*_{c,b}(G_1\times
G_2)\ar[uu]_c &&& \\
&&H^*_{c,b}(G_2).\ar@{_{(}->}[lu]_{\cong}\ar[uu]_c && } \label{Equ:
big diagram}%
\end{equation*}
Here, all vertical arrows are the natural comparison maps between
continuous bounded and continuous cohomology groups. The horizontal
maps are induced from $\varphi,\rho,\overline{\rho}$ and the restriction of $\rho$ to $\Gamma_2$.
The upwards diagonal maps are induced by the canonical projection on
the second factor.

Since the Euler class $\varepsilon_{m}\in H_{c}^{m}(\mathrm{GL}^+%
(m,\mathbb{R}))$ is bounded \cite{IvTu82}, it follows from the
commutativity of the above diagram and the isomorphism
$H_{c,b}^{m}(G_{1}\times G_{2})\cong H_{c,b}^{m}(G_{2})$ that its
image in $H_{c}^{m}(G_{1}\times G_{2})$ is in the image of the map
$H_{c}^{m}(G_{2})\hookrightarrow H_{c}^{m}(G_{1}\times G_{2})$, and
furthermore $\rho^{\ast}(\varepsilon_{m})$ is in the image of
$H^{m}(\Gamma _{2})\hookrightarrow
H^{m}(\Gamma_{1}\times\Gamma_{2})$. But since $M_{2}$ is aspherical
and $m>\mathrm{dim}(M_{2})$, we have $H^{m}(\Gamma_{2})\cong
H^{m}(M_{2})=0$, and hence also $\rho^{\ast}(\varepsilon_{m})=0$.
\end{proof}


We will also make use of the following:

\begin{lem}\label{lemma: can mod out amenable} Let $G$ be a closed subgroup of $\text{GL}^{+}(m,\mathbb{R})$, and suppose that $G=S\ltimes A$ is a decomposition of $G$ into a semidirect product where $A$ is a closed amenable normal subgroup. Denote by $p$ the projection $p:S\ltimes A\rightarrow S.$ Let $\Gamma$ be a discrete group and \[ \rho:\Gamma\longrightarrow S\ltimes A<\text{GL}^{+}(m,\mathbb{R})\] a representation. Then \[ \rho^{*}(\varepsilon_{m})=(p\circ\rho)^{*}(\varepsilon_{m})\in H^{m}(\Gamma,\mathbb{R}).\] \end{lem}

\begin{proof} Let $i:S\hookrightarrow S\ltimes A$ denote the embedding of $S$ in the semidirect product of $S$ and $A$. Since $p\circ i$ is the identity on $S$, the induced map $(p\circ i)^{*}=i^{*}\circ p^{*}$ is the identity both on $H_{c}^{*}(S)$ and on $H_{c,b}^{*}(S)$.

On the bounded continuous cohomology groups, it follows from \cite[Lemma 7.5.10]{Mo01} that the projection $p$ furthermore induces an isometric isomorphism
\[
\xymatrix{p^{*}:H_{c,b}^{m}(S)\ar[r]^{\cong\ } & H_{c,b}^{m}(S\ltimes A).}
\]
Consider the following commutative diagram: 
\[
\xymatrix{ &  & H^{m}(\Gamma)\\
H_{c}^{m}(\mathrm{GL}^{+}(m,\bbr))\ar[r] & H_{c}^{m}(S\ltimes A)\ar[ur]^{\rho^{*}} & H_{c}^{m}(S)\ar@{_{(}->}_{\
\ \ \ p^{*}}[l]\ar[u]_{\rho^{*}\circ p^{*}}\\
H_{c,b}^{m}(\mathrm{GL}^{+}(m,\bbr))\ar[u]_{c}\ar[r] & H_{c,b}^{m}(S\ltimes A)\ar[u]_{c} & H_{c,b}^{m}(S).\ar[u]_{c}\ar@{_{(}->}_{\ \ \ \cong}[l].}
\]
On the continuous cohomology groups, $p^{*}$ is not necessarily the
inverse of $i^{*}:H_{c}^{m}(S\ltimes A)\rightarrow H_{c}^{m}(S)$.
However, we get from the above diagram that $p^{*}\circ i^{*}:H_{c}^{m}(S\ltimes A)\rightarrow H_{c}^{m}(S\ltimes A)$
is the identity when restricted to the image of $H_{c,b}^{m}(S\ltimes A)\rightarrow H_{c}^{m}(S\ltimes A)$.
In particular, since the Euler class is bounded, denoting by $\varepsilon_{m}$
the restriction of the Euler class both to $S\ltimes A$ and to $S$,
we obtain $\varepsilon_{m}=p^{*}(\varepsilon_{m})$ as elements of
$H_{c}^{m}(S\ltimes A)$, and hence also $\rho^{*}(\varepsilon_{m})=(p\circ\rho)^{*}(\varepsilon_{m})$.
\end{proof}


\section{Euler class norm as a proportionality constant}\label{section: Proportionality Principles}

The purpose of this section is to prove the following:

\begin{prop}\label{prop: Proportionality Principle Simpl Vol and Euler}
Let $M$ be a closed
oriented $2n$-dimensional Riemannian manifold whose universal cover is
isometric to the Cartesian product of $n$ copies of the hyperbolic plane
$\mathcal{H}$. Then%
\[
\left\Vert M\right\Vert =\frac{\chi(M)}{(-2)^{n}\left\Vert \varepsilon_{2}%
\cup\ldots\cup\varepsilon_{2}\right\Vert _{\infty}}.
\]
\end{prop}

\medskip

For a general closed oriented Riemannian manifold $M$, the
simplicial volume and the volume are related by the Gromov--Thurston
fundamental proportionality principle (see \cite{Gr82,Th78})
\begin{equation*}\label{eq:Gromov-Thurston}
 \left\Vert M\right\Vert =\frac{\mathrm{Vol}(M)}{c(\widetilde{M})},
\end{equation*}
where $c(\widetilde{M})\in \mathbb{R}_{+}\cup\{+\infty\}$ is a
constant depending only on the Riemannian universal cover
$\widetilde{M}$ of $M$.

For locally symmetric spaces of noncompact type, it is 
shown in \cite{Bu06} that the proportionality constant satisfies
\begin{equation*}\label{eq:Michelle-vol-form}
 c(\widetilde{M})=\left\Vert \omega_{\widetilde{M}}\right\Vert_{\infty},
 \end{equation*}
where $\omega_{\widetilde{M}}\in
H_{c}^{n}($Isom$(\widetilde{M})^{\circ})$ is the image under the
Van Est isomorphism of the Isom$(\widetilde{M})^{\circ}$-invariant
volume form $\omega_{\widetilde{M}}$ on the symmetric space
$\widetilde{M}$.

When $\widetilde{M}=\mathcal{H}^n$, the Euler
characteristic and the volume of $M$ are related by Hirzebruch's
proportionality principle \cite{Hi58}
\begin{equation*}\label{eq:Hirzebruch}
 \mathrm{Vol}(M)=(-2\pi)^{n}\chi(M).
\end{equation*}
Thus, for such manifolds we have
\[
\left\Vert M\right\Vert =\frac{\mathrm{Vol}(M)}{\left\Vert\omega_{\mathcal{H}^{n}}\right\Vert _{\infty}}= \frac{(-2\pi)^n\chi(M)}{\left\Vert\omega_{\mathcal{H}^{n}}\right\Vert _{\infty}}.
\]

In order to prove the proposition, we need to understand the relation between the norms of the
volume form and of the Euler class. On the one hand, since the projection
$\prod_{i=1}^{n}\mathrm{SL}(2,\mathbb{R})\rightarrow
\prod_{i=1}^{n}\mathrm{PSL}(2,\mathbb{R})$ is finite-to-one, the
volume form
$\omega_{\mathcal{H}^{n}}$ defines
continuous cohomology classes with the same sup norm in
$H^{2n}_c(\prod_{i=1}^{n}\mathrm{PSL}(2,\mathbb{R}))$ and in
$H^{2n}_c(\prod_{i=1}^{n}\mathrm{SL}(2,\mathbb{R}))$. 
On the other hand, since $\GL^+(2,\bbr)\cong \SL(2,\bbr)\times \bbr^{>0}$, by
Lemma \ref{lemma: can mod out amenable} the sup norm of the Euler class
is unchanged by replacing it with its restriction to $\prod_{i=1}^n\SL(2,\BR)$ which we still denote by $\varepsilon_2\cup\ldots\cup \varepsilon_2$.
The proposition will follow from:

\begin{lem}\label{lem:Euler-volume} 
The volume form and the Euler class are related by
$$ 
\omega_{\mathcal{H}^n}=(-4\pi)^n \varepsilon_2\cup\ldots\cup \varepsilon_2
$$
as elements in $H^{2n}_c(\prod_{i=1}^{n}\mathrm{SL}(2,\mathbb{R}))$.
\end{lem}

\begin{proof}[Proof of Lemma \ref{lem:Euler-volume}]
As the Van Est isomorphism is
multiplicative, we have
$$
 \omega_{\mathcal{H}\times\ldots\times\mathcal{H}}=\omega_{\mathcal{H}}\cup\ldots\cup\omega_{\mathcal{H}}\in H_{c}^{2n}
 (\prod_{i=1}^{n}\mathrm{SL}(2,\mathbb{R})),
$$
so it is enough to show that
\begin{equation}\label{eq:surface-vol-Euler}
\omega_{\mathcal{H}}=-4\pi \varepsilon_2
\end{equation}
as elements of $H_c^2(\SL(2,\bbr))$.

Let $\gS_g$ be a closed surface of genus $g\geq 2$ and let $\rho:\pi_{1}(\Sigma_{g})\rightarrow\mathrm{PSL}%
(2,\mathbb{R})$ be the embedding of $\pi_{1}(\Sigma_{g})$ corresponding to an arbitrary complete hyperbolic structure of $\Sigma_{g}$. It is shown in \cite{Mi58} (see also \cite[pages 312-314]{MiSt79} for an alternative
topological proof) that $\rho$ lifts to a
representation $\overline{\rho}:\pi_{1}(\Sigma_{g})\rightarrow
\mathrm{SL}(2,\mathbb{R})$ and that
$$
 \left\langle\overline{\rho}^{\ast }(\varepsilon_{2}),[\Sigma_{g}]\right\rangle
 =1-g=\frac{\chi(\Sigma_g)}{2}.
$$
Since $H^{2}(\pi_{1}(\Sigma_{g}))$ is one dimensional, and
$$
 \left\langle\overline{\rho}^{\ast }(\omega_{\mathcal{H}}),[\Sigma_{g}]\right\rangle=
 \left\langle\rho^{\ast}(\omega_{\mathcal{H}}),[\Sigma_{g}]\right\rangle =\mathrm{Vol}
 (\Sigma_{g})=-2\pi\chi(\Sigma_{g}),
$$
we have $\overline{\rho}^{\ast}(\omega_{\mathcal{H}})=-4\pi
\overline{\rho}^{\ast}(\varepsilon_{2})$. Moreover, as $\overline{\rho}(\gS_g)$ is a cocompact lattice of $\SL(2,\BR)$, 
$\overline{\rho}^{\ast}$ induces an isomorphism
between $H_{c}^{2}(\mathrm{SL}(2,\mathbb{R}))$ and $H^{2}(\pi_{1}%
(\Sigma_{g}))$, and hence
$\omega_{\mathcal{H}}=-4\pi \varepsilon_{2}$. This completes the
proof of Lemma \ref{lem:Euler-volume}, and hence of Proposition
\ref{prop: Proportionality Principle Simpl Vol and Euler}.
\end{proof}

\begin{rem}
By Proposition \ref{prop: Proportionality Principle Simpl Vol and Euler}, 
computing the simplicial volume of a closed $\mathcal{H}^n$-manifold is equivalent to
computing the sup norm of $\varepsilon_2\cup\ldots\cup
\varepsilon_2\in H^{2n}_c(\Pi_{i=1}^n \GL^+(2,\bbr))$. Milnor's
original inequality \cite{Mi58} amounts to showing
$\left\Vert\varepsilon_2\right\Vert_\infty=1/4$, and in \cite{Bu07},
the first author proved $\left\Vert \varepsilon_2 \cup
\varepsilon_2\right\Vert_\infty = 1/24$. For more than two factors,
only rough lower and upper bounds are currently known.
\end{rem}


\section{Representations of product groups}\label{Section: Representations of
products}

In this section we prove the following general proposition\footnote{Unfortunately the variant of this proposition that we included in the announcement \cite[Lemma 3.1] {BuGe08} is erroneous.}

\begin{prop}\label{lem:prod}\label{prop:prod}
Let $H=\prod_{i=1}^n H_{i}$ be a direct product of $n$ groups and let $\rho:H\to
\text{GL}^{+}(m,{\mathbb{R}})$ be an oriented representation such that
$\rho(H_{i})$ is nonamenable for each $i$. Then $m\geq2n$. When
$m=2n$, 
the identity component of the Zariski closure of $\rho(H)$ is reductive, and $\rho$ factors through a map to $\prod_{i=1}^n\GL(2,\BR)$
\begin{equation*}
\xymatrix{ \prod_{i=1}^{n}H_{i}
\ar@/^1.5pc/[rr]^{\rho}\ar[r]^{{\rho_1}} & \prod_{i=1}^n\GL(2,\BR)
\ar[r]^{\rho_2} & \GL^+(2n,\bbr),}
\end{equation*}
and there is an even integer $0\le t\le n$ such that, up to reordering the factors $H_i$, $\rho_1=\prod_{i=1}^n\rho_{1,i}$ where $\rho_{1,i}$ is a representation of $H_i$ and the restriction of its image in $\prod_{i=1}^n\GL_2(\BR)$ to the $j$'th factor is irreducible when $j=i$ and scalar otherwise, and $V={\mathbb{R}}^{2n}$ decomposes as an invariant direct sum 
$$
 V=\bigoplus_{i=1}^{t/2}W_i\oplus\bigoplus_{i=t+1}^nV_i,
$$
where the $W_i$'s are $4$-dimensional and the $V_i$'s are $2$-dimensional. 
The restriction of $\rho_2$ to any pair of factors corresponding to indices $(2i-1,2i)$ with $i\le t/2$ on the invariant subspace $W_i$ is isomorphic to the tensor of two standard representations of $\GL_2(\BR)$ on $\BR^2$,  and is 
trivial on any of the other $W_j$ and $V_j$. The restriction of $\rho_2$ to a factor corresponding to $i>t$ and to $V_i$ is isomorphic to the standard representation of $\GL(2,\BR)$ on $\BR^2$, while it is 
trivial on any other $V_j$ and $W_j$. The decomposition is uniquely determined, and in particular the number $t=t(\rho)$ is well defined. 
When $t=0$, $\rho(H)$ is conjugated to a subgroup of $\text{Bl}_{2,n}^+$.\end{prop}

The following classical generalization of Schur's lemma is implicit in the proof of Theorem 3 in \cite{Clifford} when restricting to the case of product groups:

\begin{lem}\label{lem:Clifford}
Let $r:H_1\times H_2\to\GL(W)$ be an irreducible representation of a direct product over an algebraically closed field. Then $r$ is isomorphic to the tensor of some (irreducible) representations $r_i:H_i\to \GL(\ti W_i)$, $i=1,2$.  

Moreover, any linear automorphism $T\in GL(W)$ which commutes with $r(H_2)$ is of the form $\ti T\otimes 1$ for some $\ti T\in\GL (\ti W_1)$. 
\end{lem}

\begin{proof}
Let $U\subset W$ be an $r(H_1)$ nontrivial irreducible subspace. If $U=W$, the lemma follows directly from Schur's lemma with $\ti W_1 =W$ and $\dim(\ti W_2)=1$. Suppose that $U\ne W$. Since $W$ is irreducible, we can find $h_1,\ldots,h_k\in H_2$ with $h_1=1$ such that $W=\sum_{i=1}^kr(h_i)U$ and $r(h_j)U$ is not contained in $\sum_{i=1}^{j-1}r(h_i)U$ for every $j\le k$. Since the $h_i$ commute with $H_1$, $r(h_j)U$ are $H_1$ isomorphic invariant subspaces, hence $H_1$-irreducible, so we conclude that $W=\bigoplus_{i=1}^k r(h_i)U$ is a direct sum. Choose an ordered basis $B_1$ of $U$ and extend it to an ordered basis $B$ of $W$ by adding $r(h_i)B_1,~i=2,\ldots,k,$ to it. 
With respect to the basis $B$ the elements of $r(H_1)$ are block diagonal with $k$ identical blocks. To prove the lemma, we have to show that if $T$ is a linear transformation of $W$ commuting with $r(H_1)$ then each of the $k\times k$ blocks of $T$ determined by the basis $B=B_1\cup\ldots \cup B_k$ is a scalar matrix. In other words, denoting by $P_j:W\to r(h_j)U,~j=1,\ldots,k$, the associated projections, we have to show that $r(h_j^{-1})\circ P_j\circ T\circ r(h_i)$ is a scalar transformation of $U$ for every $i,j$. Since all four transformations involved commute with $r(H_1)$, and $U$ is $H_1$ irreducible, this is a consequence of Schur's lemma. 
\end{proof}

The following simple observation will be used in the proof below. 

\begin{lem}\label{lem:unique-decomposition}
Suppose that $\gD$ is a subgroup of $GL(m,\BC)$ and $\BC^m$ decomposes as a direct sum of one $\gD$-invariant subspace $U$ of dimension $>1$ and $(m-\dim U)$ $\gD$-invariant lines. Then every $\gD$-irreducible subspace of dimension $>1$ is contained in $U$.
\end{lem}

\begin{proof}
Let $W$ be a $\gD$-invariant subspace of dimension $>1$. If $W$ is not contained in $U$, then the corresponding projection $P_L$ on one of the lines $L$ in the given decomposition is nonzero on $W$. This implies that $W$ is not irreducible since $\ker (P_L|_W)$ is a $\gD$-invariant subspace of $W$ of codimension $1$.  
\end{proof}

\begin{proof}[Proof of Proposition \ref{prop:prod}]
Consider $V_{\BC}=\BC\otimes V$ as a complex vector spaces. By choosing a basis for $V$ we may identify 
$\text{GL}(m,{\mathbb{C}})$ with the group of complex linear automorphisms of $V_\BC$ and $\GL(m,\BR)$ as the subgroup preserving the real subspace $V$.

Let $G_{i}$ be the Zariski closure of $\rho(H_{i})$ in $\text{GL}%
(m,{\mathbb{C}})$. 
Let $G_i^\circ$ be the identity connected component of $G_i$ with respect to the Zariski topology, and let
$G_{i}^\circ=S_{i} \ltimes U_{i}$ be a Levi decomposition of the complex algebraic group $G_{i}^\circ$, $S_{i}$ being reductive and $U_{i}$ unipotent (see \cite{Margulis}, 0.28).
The $G_{i}$ commute with each other, being the Zariski closure of commuting groups. Hence the $S_i$ also commute with each other. Let $S=\prod_{i=1}^{n}S_{i}$ be their product. Then $S$ is reductive. Our nonamenability assumption implies in particular that
each $S_{i}$ is nonabelian.

We will first prove the analog of the statement when $\prod H_i$ is replaced by $S=\prod S_i$, $\rho$ by the {\it complex} representation corresponding to the inclusion of $S$ in $\GL(m,\mathbb{C})$, and the nonamenability assumption by the data that each $S_i$ is nonabelian. That is, we will show that $m\ge 2n$ and in the equality case, $V_\BC$ decomposes as an $S$-invariant direct sum $\bigoplus_{i=1}^{t/2}W_{\BC,i}\oplus\bigoplus_{i=t+1}^nV_{\BC,i}$ where the $W_{\BC,i}$ are $4$-dimensional and isomorphic to tensors of pairs of invariant complex planes corresponding to pairs of factors of $S$, the $V_{\BC,i}$ are planes, and each $S_i$ acts nonscalarly on exactly one subspace in this decomposition. Moreover we will show that this decomposition is unique. 

To prove this claim we argue by induction on $n$, the case $n=1$ being trivial. Suppose $n>1$.
Since $S$ is reductive, $V_{\BC}$ decomposes as an $S$-invariant direct sum of $S$-irreducible
subspaces $\bigoplus_{j=1}^{l}\mathcal{V}_{{\mathbb{C}},j}$. 

Suppose first that $l>1$. Then since each $S_i$ is nonabelian, its restriction to at least one of the $\mathcal{V}_{{\mathbb{C}},j}$ is nonabelian, and the conclusion that $m\ge 2n$ follows from the induction hypothesis applied to all the $\mathcal{V}_{{\mathbb{C}},j}$ simultaneously, i.e.
$$
 m=\sum _{j=1}^{l}\dim \mathcal{V}_{{\mathbb{C}},j}\ge 2\sum _{j=1}^{l} \#\{i:\text{the restriction of}~S_i~\text{to}~\mathcal{V}_{{\mathbb{C}},j}~\text{is nonabelian}\}\ge 2n.
$$
Moreover if $m=2n$, each $S_i$ is nonabelian on exactly one $\mathcal{V}_{{\mathbb{C}},j}$, and since the restriction of $S_i$ to every other $\mathcal{V}_{{\mathbb{C}},j}$ (i.e. one on which it acts commutatively) is central, it follows from Schur's lemma that it acts scalarly there. The uniqueness of the decomposition $V_\BC=\bigoplus_{i=1}^{t/2}W_{\BC,i}\oplus\bigoplus_{i=t+1}^nV_{\BC,i}$ follows in this case from the induction hypothesis and Lemma \ref{lem:unique-decomposition}. 

Suppose now that $l=1$, i.e. $S$ acts irreducibly on $V_\BC$. Writing $S=S_1\times\prod_{i=2}^nS_i$, we derive from Lemma \ref{lem:Clifford} that $V_\BC\cong \ti{V}_1\otimes \ti{V}_2$ where $\ti{V}_1$ is a representation of $S_1$ and $\ti{V}_2$ is a representation of $\prod_{i=2}^nS_i$. The assumption that all $S_i$ are nonabelian and the inductive hypothesis give $\dim \ti{V}_1\ge 2$ and $\dim \ti{V}_2\ge 2(n-1)$. Thus $m=\dim V_\BC\ge 4(n-1)\ge 2n$, and equality implies that $n=2,m=4$ and our representation is isomorphic to the tensor of two $2$-dimensional representations of $S_1$ and $S_{2=n}$. This finishes the proof of the analogous statement for the complex representation given by the inclusion of $S$ in $\GL(m,\mathbb{C})$.

We shall denote by $j^\circ$ and $k^\circ$ either a pair of indices $(2i-1,2i)$ or a single index $i$ according to whether $i\le t/2$ or $i>t$. We let $|j^\circ|=1,2$ according to whether $j^\circ$ is a single index or a pair of indices.
For each pair $j^\circ=(2i-1,2i),~1\le i\le t/2$, we let 
$$
 G^\circ_{j^\circ}=G^\circ_{2i-1}\times G^\circ_{2i},~S_{j^\circ}=(S_{2i-1}\times S_{2i})~\text{~and~}~U_{j^\circ}=(U_{2i-1}\times U_{2i}),
$$
and let $V_{\BC,j^\circ}$ be the unique $S_{j^\circ}$-irreducible space of dimension $>1$ (Lemma \ref{lem:unique-decomposition}).

Consider now a unipotent element $u\in U_{k^\circ}$. 
Since for $j^\circ\ne k^\circ$, $u$ centralizes
$S_{j^\circ}$, it preserves $V_{\BC,j^\circ}$, and by Schur's lemma acts scalarly, hence, being unipotent, trivially on
$V_{{\mathbb{C}},j^\circ}$.
To show that $V_{{\mathbb{C}},k^\circ}$ is $u$
invariant as well, we argue by way of contradiction. Suppose that $u(v)\notin
V_{{\mathbb{C}},k^\circ}$ for some $v\in V_{{\mathbb{C}},k^\circ}$. Let $j^\circ\ne k^\circ$ be such that the projection of $u(v)$ to $V_{\BC,j^\circ}$ is nontrivial, and choose $s\in S_{j^\circ}$ for which $su(v)$ is not on the
line spanned by $u(v)$. Such $s$ exists since $S_{j^\circ}$ acts irreducibly on $V_{\BC,j^\circ}$.
This however gives the impossible
\[
\begin{array}{lll}
2 &=& \text{rank}\{su(v),u(v)\} = \text{rank}\{ u^{-1}su(v),v\}\\
   &=& \text{rank}%
\{s^{-1}u^{-1}su(v),s^{-1}(v)\}=1,
\end{array}
\]
where the last equality holds since $s$ and $u$ commute and $s$ acts scalarly on
${\mathbb{C}} v$. This holds for any $u\in U_{k^\circ}$, so all the $V_{\BC,j^\circ}$ are $G_{k^\circ}=S_{k^\circ} \ltimes U_{k^\circ}$-invariant, and since $k^\circ$ is arbitrary, it follows that the decomposition $V_\BC=\bigoplus V_{\BC,j^\circ}$ is $G=S\ltimes U$ invariant.

Next we show that $G^\circ$ is reductive.
Consider first $k^\circ$ with $|k^\circ|=1$. Since $U_{k^\circ}$ is unipotent, the subspace of $U_{k^\circ}$-invariant vectors in $V_{\BC,k^\circ}$ is non-trivial. As $U_{k^\circ}$ is normalized by $S_{k^\circ}$, this subspace is $S_{k^\circ}$-invariant, and as $S_{k^\circ}$ is irreducible on $V_{\BC,k^\circ}$, it must be the full space, i.e. $U_{k^\circ}$ acts trivially on $V_{\BC,k^\circ}$. We already saw in the previous paragraph that $U_{k^\circ}$ acts trivially on $V_{\BC,j^\circ}$ for any $j^\circ\ne k^\circ$, so $U_{k^\circ}$ is trivial. 
Similarly, if $|k^\circ|=2$, we write $G^\circ_{k^\circ}=G^\circ_{2i-1}\times G^\circ_{2i}=(S_{2i-1}\times S_{2i})\ltimes (U_{2i-1}\times U_{2i})$ and $V_{\BC,k^\circ}=\ti V_{2i-1}\otimes \ti V_{2i}$ where the $\ti V_j$ are irreducible $G_j$ spaces, $j\in\{2i-1,2i\}$. If $u$ is a unipotent belonging to one of the factors, say $u\in U_{2i-1}$, then, being in the centralizer of the other factor, $u$ preserves the tensor structure of $V_{\BC,k^\circ}$ and acts trivially on its second factor (Lemma \ref{lem:Clifford}). As above we deduce that the space of $U_{2i-1}$ invariants in $\ti V_{2i-1}$ is the full space, and hence that $U_{2i-1}$ acts trivially on $V_{\BC,k^\circ}$ as well as on $V_{\BC,j^\circ}$ for any $j^\circ$.
Thus $U$ is trivial, i.e. $G^\circ=S$ is reductive.

Finally, since the decomposition established above is unique, and the $G^\circ_{k^\circ}$ are defined over ${\mathbb{R}}$, each $V_{{\mathbb{C}},k^\circ}$ is invariant under complex conjugation, hence is defined over ${\mathbb{R}}$. In other words $V_{{\mathbb{C}},k^\circ}=\BC\otimes V_{k^\circ}$ -- the complexification of the real space $V_{k^\circ}=V_{{\mathbb{C}},k^\circ}\cap V$. The uniqueness of the decomposition also implies that the result holds for $G=\prod G_i$ rather than its Zariski identity component $\prod G_i^\circ$.

Moreover, in the case $t=0$, since the representation $\rho$ is oriented and the elements of $\rho(H_i)$ act scalarly and hence preserve orientation on the real $2$-dimensional spaces $V_j$, for $j\neq i$, they also preserve the orientation on $V_i$. This proves the last statement of the proposition.
\end{proof}

\begin{rem}
$(i)$ In later sections we shall need to apply Proposition \ref{prop:prod} in a slightly more general setup where $H$ is an almost direct product of $n$ factors, i.e $H=(\prod_{i=1}^nH_i)/C$ where $C\lhd \prod_{i=1}^nH_i$ is a (finite) central subgroup. In this case, given $\rho:H\to\GL^+(2n,\BR)$ and letting $\ti\rho:\prod_{i=1}^nH_i\to\GL^+(2n,\BR)$ be the induced representation, we can decompose $\ti\rho=\ti\rho_2\circ\ti\rho_1$, as in Proposition \ref{prop:prod}. Then $\ti \rho_1(C)$ lies both in the center of $\prod_{i=1}^n\GL(2,\BR)$ and in the kernel of $\ti\rho_2$, and $\rho$ decomposes as $\rho=\rho_2\circ\rho_1$ through a map to $\prod_{i=1}^n\GL(2,\BR)/\rho_1(C)$ which is an almost direct product. If the restriction of $\ti\rho_2$ to a factor or a pair of factors of $\prod_{i=1}^n\GL(2,\BR)$ is irreducible or tensor of irreducible or scalar on some subspace of $V=\BR^{2n}$ the same holds for the restriction of $\rho_2$ to the image of that factor (or pair of factors) in $\prod_{i=1}^n\GL(2,\BR)/\rho_2(C)$. We shall denote $t(\rho)=t(\ti\rho)$. 

$(ii)$ In the case when $S\le\GL^+(2n,\BR)$ is an almost direct product of $2n$ nonamenable subgroups, we shall simply write $t(S)=t(i)$, where $i:S\hookrightarrow \GL^+(2n,\BR)$ is the inclusion.
\end{rem}

We end this section with an extension of Lemma \ref{lemma: tensor product} that will be applied repeatedly in the sequel combined with Proposition \ref{prop:prod}.

\begin{lem} \label{lemma: tensor product factor}\label{lem:trho=0}
Let $\rho: \prod_{i=1}^n H_i \rightarrow \GL^+(2n,\mathbb{R})$ be a representation such that $\rho(H_i)$ is nonamenable for every $i$. If $t(\rho)>0$ then
$$ 
 \rho^*(\varepsilon_{2n})=0.
$$
\end{lem}

\begin{proof} 
By Proposition \ref{prop:prod}, $\rho$ factors through a map to $\prod_{i=1}^n\GL_2(\BR)$:
\begin{equation*}
\xymatrix{ \prod_{i=1}^{n}H_{i}
\ar@/^1.5pc/[rr]^{\rho}\ar[r]^{{\rho_1}} & \prod_{i=1}^n\GL_2(\BR)
\ar[r]^{\rho_2} & \GL^+(2n,\bbr).}
\end{equation*}
Up to replacing $\prod_{i=1}^n H_i$ by a finite index subgroup, which has  no influence on the vanishing of the Euler class, we can assume that the image of $\rho_1$ is contained in the product  $\prod_{i=1}^n\GL^+_2(\BR)$. If $t(\rho)>0$, then up to permuting the factors, Proposition \ref{prop:prod} allows us further to describe $\rho_2$ as follows: 
Let $\iota$ denote the canonical inclusion of $\prod_{i=1}^n\SL_2(\BR)$ in $\prod_{i=1}^n\GL^+_2(\BR)$. 
Precomposed with $\iota$, the representation $\rho_2$ has the form 
$$
 \rho_2\circ \iota (A,B,C)= \left( \begin{array}{cc} A\otimes B & 0 \\ 0 & \rho_2'(C) \end{array} \right),
$$
for $A,B\in \SL_2(\mathbb{R})$ and $C\in\Pi_{i=3}^n\SL_2(\BR)$, where
$\rho_2':\Pi_{i=3}^n\SL_2(\BR)\to\GL^+_{2n-4}(\mathbb{R})$ is some representation. Denote by $p:\Pi_{i=1}^n\GL^+_2(\BR)\rightarrow \Pi_{i=1}^n\SL_2(\BR)$ the natural projection. Note that $\rho_2$ and $\rho_2 \circ \iota \circ p$ lie in the same path connected component in the space of representations, thus $\rho_2^*(\varepsilon_m)=(\rho_2 \circ \iota \circ p)^*(\varepsilon_m)$. 
Finally, since $\rho_2\circ \iota$ is the inclusion of the direct sum of the tensor representation $\rho_\otimes$ of $\SL(2,\mathbb{R})\times \SL(2,\mathbb{R})$ and the representation $\rho_2'$ on $\Pi_{i=3}^n\SL_2(\BR)$, we have 
$$(\rho_2 \circ \iota)^*(\varepsilon_m)=\rho_\otimes^*(\varepsilon_4 )\cup (\rho_2')^*(\varepsilon_{m-4})=0,$$
since $\rho_\otimes^*(\varepsilon_4 )=0$ by Lemma \ref{lemma: tensor product}.
\end{proof}




\section{Lattices in $\PSL(2,\mathbb{R})^n$ and their representations}\label{Section: Representations}

Recall that a lattice in a semisimple Lie group $G$ is said to be {\it irreducible} if it projects densely to every proper quotient of $G$.
Let $G=\prod_{i=1}^{n}\text{PSL}(2,\bbr)\cong\text{Isom}(\mathcal{H}^n)^\circ$, and let $\Gamma$ be a lattice\footnote{In this section we do not need to assume that $\gC$ is cocompact.} in $G$.
By \cite[Theorem 5.22]{Raghunathan}, up to replacing $\gC$ by a finite index subgroup $\gC'\le \gC$, $G$ decomposes as a direct product $G=\prod_{j=1}^mG_j$ such that:
\begin{itemize}
\item each $G_j$ is a direct product of $\PSL(2,\BR)$'s,
\item $\gC_j':=\gC'\cap G_j$ is an irreducible lattice in $G_j$, and
\item $\gC'$ is the direct product of the $\gC_j'$'s.
\end{itemize}
For the rest of this section let us assume that $\gC=\gC'$, i.e. that $\gC$ itself decomposes as a direct product of irreducible factors $\gC_j=\gC\cap G_j$.
Let $\rho:\Gamma\to\text{GL}^{+}(2n,\bbr)$ be an orientable
$2n$-dimensional representation of $\Gamma$. In order to analyze the possible images of $\rho$, we shall distinguish between 3 cases:

\medskip

\noindent\textbf{Case 1: $\Gamma$ is completely reducible}, i.e. $\Gamma
=\prod_{i=1}^{n} \Gamma_{i}$, where each $\gC_i$ is a lattice in $\PSL(2,\BR)$.
It is an immediate consequence of Proposition \ref{prop:prod} that:

\begin{prop}\label{prop: Case 1 Gamma completely reducible}
In this case either
\begin{itemize}
\item $\rho(\gC_i)$ is amenable for some $i=1,\ldots,n$, 
\item $t(\rho)>0$, or
\item $\rho (\gC)$ is conjugate to a subgroup of $\text{Bl}_{2,n}^+$.
\end{itemize}
\end{prop}

\medskip


\noindent\textbf{Case 2: $\Gamma$ is rigid}, i.e. each $G_j$ has real rank at least $2$, or equivalently no $G_j$ is $\PSL(2,\BR)$.
By Margulis arithmeticity theorem, every $\gC_i$ is a Hilbert--Blumenthal modular group. 

\begin{prop}\label{prop:rigid}
In this case, up to replacing $\gC$ by a finite index subgroup, $\rho(\Gamma)$ is contained in a connected semisimple Lie group $S$ for which all the noncompact simple factors are locally
isomorphic to $\text{PSL}(2,\mathbb{R})$.
The number of
the simple factors of $S$ is at most $n$, and in the case it is
exactly $n$, either 

\noindent $(a)$ $t(S)>0$, or

\noindent $(b)$ $S$ has no compact factors, and

\begin{enumerate}
\item[($b_1$)] it is conjugate to $\text{SBl}_{2,n}$, the diagonal product of $\SL_2(\BR)'s$,
\item[($b_2$)] $\rho$ is faithful and $\rho(\gC)\cap Z(S)=\emptyset$, and
\item[($b_3$)] the image of $\rho(\Gamma)$ in $S$ is a lattice and is cocompact in the case $\gC\leq G$ is cocompact.
\end{enumerate}
\end{prop}

\begin{proof}
By Whitney's theorem (see \cite[Theorem 3.6]{Pl-Ra}), the Zariski closure
of $\rho(\Gamma_{j})$ in $\text{GL}^+(2n,{\mathbb{R}})$, which we shall denote below by $\overline{\rho(\Gamma_{j})}^z$, has finitely many connected components with respect to the Hausdorff topology induced from $\text{GL}^+(2n,{\mathbb{R}})$. Hence,
up to replacing $\Gamma$ by a finite index subgroup, we may assume that 
$\overline{\rho(\Gamma_{j})}^z$ is connected (as a real group) for each $j\le m$. Moreover, as $\overline{\rho(\Gamma_{j})}^z$ is closed in $\GL^+(2n,\BR)$, it is a Lie group (see \cite[Proposition 1.75]{Knapp}).


Since for each $j$, $\rank_\BR(G_j)\ge2$, it follows from Margulis' Theorem (see \cite[Ch. IX, Theorem 5.8]{Margulis}) that $\overline{\rho(\gC_j)}^z$ is semisimple, and hence $\overline{\rho(\gC)}^z$ is also semisimple. Lets denote $S=\overline{\rho(\gC)}^z$, and let $S_{nc},S_c$ be the product of noncompact, resp. compact, simple factors of $S$. Then $S$
is an almost direct product $S=S_{nc}\times S_c$ (i.e. $S_{nc}\cap S_c$ is finite and central in $S$).

Let $\tilde{\rho}:\Gamma_{j}\to
\text{Ad}(S_{nc})$ be the representation induced from $\rho$ by dividing out $S_c$ and composing with the Adjoint representation of $S_{nc}$.
Since each simple factor $S^{\prime
}$ of $S_{nc}$ is noncompact, the projection of the image of each $\ti\rho
(\Gamma_{j})$ to $S^{\prime}$, being normal, is either trivial or Zariski dense,
and $\prod_{j=1}^m\Gamma_{j}$ is superrigid in
$G=\prod_{j=1}^mG_j=\prod_{i=1}^n\text{PSL}(2,{\mathbb{R}})$ (as $\rank_\BR(G_j)\ge2,~\forall j$), it
follows from Margulis' Superrigidity Theorem \cite[Ch. VII, Sec. 5]{Margulis}
that $\ti\rho$ extends to a representation of $G=\prod_{i=1}^n\text{PSL}(2,{\mathbb{R}})$. Thus $\text{Ad}(S_{nc})$ is a homomorphic image of $\prod_{i=1}^n\PSL(2,\BR)$. This proves the first statement of the proposition.
 
In order to see that the number of simple factors of $S$ is at most $n$, we apply Proposition \ref{prop:prod} to the product of the simple factors of $\ti S$, the universal covering Lie group of $S$. Note that all the simple factors of $S$ (including the compact ones) considered as abstract groups with no topology are nonamenable. Suppose now that $S$ has exactly $n$ simple factors. Then if $t(S)>0$ we are in case $(a)$. Assuming $t(S)=0$, we get from Proposition \ref{prop:prod} that $S$ is conjugate to a subgroup of $\text{Bl}_{2,n}^+$.
Then all the simple factors of $S$ are locally isomorphic to $\PSL_2(\BR)$ (i.e. $S=S_{nc}$), and hence $\dim(S)=3n$. Furthermore, since $S$ is connected and semisimple, its conjugate in $\text{Bl}_{2,n}^+$ is contained in, and hence, by dimension equality, coincides with the commutator group of $\text{Bl}_{2,n}^+$ which is $\text{SBl}(2,\BR)$. This proves $(b_1)$.

$(b_2)$ follows from Margulis' Superrigidity Theorem as the map $G\to\text{Ad}(S)$ extending the map $\gC\ni\gc\mapsto \text{Ad}_S(\rho(\gc))$ is an isomorphism.

The fact $(b_3)$, that the image of $\gC$ in $S$ is a (cocompact) lattice, holds because $G/\gC$ is (equivariantly) homeomorphic to $\text{Ad}(S)/\text{Ad}_S(\rho(\gC))$ and $S$ has a finite center by $(b_1)$.
\end{proof}

\medskip


\noindent\textbf{Case 3: The mixed case}.
The mixed case is when some, but not all, of the $G_j$ are isomorphic to $\PSL(2,\mathbb{R})$. Denote by $\mathcal{R}\subset\{1,\ldots,m\}$ the set of indices corresponding to rigid factors of $\gC=\prod_{j=1}^m\gC_j$.
In this case we have:

\begin{prop}\label{prop:mixed}
Up to replacing each rigid factor by a finite index subgroup (without changing the nonrigid factors), we have that either

\noindent $(i)$ for some $j\notin\mathcal{R}$, $\rho(\Gamma_{j})$ is amenable, 

\noindent $(ii)$ for some $j\in\mathcal{R}$, $\rho(\Gamma_{j})$ is contained in a group locally isomorphic
to $\prod_{i=1}^{k}\text{PSL}(2,{\mathbb{R}})$ where $k$ is strictly
smaller than the number of factors of $G_{j}$, or
  
\noindent $(iii)$ there are $n$ groups $H_i,~i=1,\ldots,n$ and an almost direct product $H=(\prod_{i=1}^nH_i)/C$ such that $\rho$ factors through a map 
\begin{equation*}
\xymatrix{ \gC
\ar@/^1.5pc/[rr]^{\rho}\ar[r]^{{\rho_1}} & H
\ar[r]^{\rho_2} &\GL^+(2n,\bbr),}
\end{equation*}
and we have that either
\begin{enumerate}
\item[($iii'$)] $t(\rho_2)>0$, or
\item[($iii''$)] $\rho(\gC)$ is conjugated to a subgroup of $Bl_{2,n}^+$.
\end{enumerate}
\end{prop}




\begin{proof}
Fix a labeling of the simple factors of $\PSL(2,\BR)^n$ by $\{1,\ldots,n\}$, and for $j\in\mathcal{R}$ denote by $K_j$ the set of indices corresponding to the factors of $G_j$ under this labeling.
By replacing every $\gC_j,~j\in\mathcal{R}$, by a finite index subgroup, we may assume that the Zariski closure $F_j=\overline{\rho(\gC_j)}^z$ of $\rho(\Gamma_{j})$ in
$\text{GL}(2n,{\mathbb{R}})$ is connected in the Hausdorff topology.
Then, as in Case (2), for each $j\in\mathcal{R}$, $F_j$ is locally isomorphic to a product of $\PSL(2,\BR)$'s, and the number of factors is at most $\dim(G_j)/3$. 
Assuming that $(ii)$ above is not satisfied, for every rigid component $j$, $F_j$ is an almost direct product of exactly $|K_j|$ such factors, so let us name its simple factors by $H_i,~i\in K_j$, and write $F_j=\prod_{i\in K_j}H_i/C_j$ where $C_j$ is a finite central group in $\prod_{i\in K_j}H_i$. Set $C=\prod_{j\in\mathcal{R}}C_j$, for $j\notin\mathcal{R}$, set $H_i=F_j=\rho(\gC_j)$, where $i=i(j)\in\{1,\ldots,n\}$ is the labeling index of $G_j$, and define $H=\prod_{i=1}^nH_i/C=\prod_{j=1}^mF_j$. Since the $\gC_j$'s commute which one another, the $F_j$'s commute and the product of the inclusion maps $F_j\to\GL^+(2n,\BR)$ gives a well defined representation, which we denote by $\rho_2:H\to\GL^+(2n,\BR)$. Obviously we set $\rho_1:\gC\to \prod_{j=1}^mF_j=H$ to be the product of the restrictions $\rho|_{\gC_j}:\gC_j\to F_j$. Finally if we assume in addition that $(i)$ is not satisfied, we may apply Proposition \ref{prop:prod} and conclude that either $(iii')$ or $(iii'')$ must hold.
\end{proof}

\section{An inequality of Euler classes}\label{Section:Estimating the norm of the Euler class} In this
section we prove the following:

\begin{thm}
\label{Thm: norm of the Euler class}
Let $\Gamma$ be a cocompact lattice in 
$\text{Isom}(\mathcal{H}^n)^+$ and
\[
\rho:\Gamma\longrightarrow\mathrm{GL}^+(2n,\mathbb{R})
\]
a representation. Then
\[
\left\Vert \rho^{\ast}(\varepsilon_{2n})\right\Vert _{\infty}\leq
\|\underbrace{\varepsilon_{2}\cup\ldots\cup\varepsilon_{2}}_{n\text{ times}%
}
\|_\infty.
\]
\end{thm}

\begin{proof}
The inclusion of a finite index subgroup $\Delta$ in $\Gamma$
induces isometric embeddings
$$H^*(\Gamma)\longrightarrow H^*(\Delta) \text{ \ \ and \ \ }H^*_{b}(\Gamma)\longrightarrow H^*_{b}(\Delta)$$
both on the standard and bounded cohomology groups. The statement
for bounded cohomology groups is a particular case of
\cite[Proposition 8.6.2]{Mo01} and the standard case is proven
identically. As a result, we may replace $\gC$ by a finite index
subgroup which is contained in the identity connected component $\text{Isom}(\mathcal{H}^n)^\circ\cong\prod_{i=1}^n\PSL(2,\BR)$ and
decomposes as in Section \ref{Section:
Representations}. We will argue case by case, showing that if
$\rho^*(\varepsilon_{2n})\ne 0$ then, up to replacing $\gC$ again by a finite index subgroup, $\rho(\gC)$
is conjugate to a subgroup of $\text{Bl}_{2,n}^+$. By Lemma
\ref{lemma: rho(Gamma) amenable}, if for some $j$, $\rho(\Gamma_j)$
is amenable, then $\rho^*(\varepsilon_{2n})=0$. We shall assume
below that this is not the case.

\noindent{\it Case 1: $\Gamma$ is completely reducible.}

In this case, assuming $\rho^*(\varepsilon_{2n})\ne 0$ we have by Lemma \ref{lem:trho=0} that $t(\rho)=0$, hence
 Proposition \ref{prop: Case 1 Gamma completely reducible}
gives that $\rho(\Gamma)$ is conjugate to a subgroup of
$\text{Bl}_{2,n}^+$.

\medskip

\noindent{\it Case 2: $\Gamma$ is rigid.}

Replacing $\Gamma$
by a further finite index subgroup, if necessary, we get from Proposition \ref{prop:rigid} that
$\rho(\Gamma)$ is contained in
a connected Lie group $S$ locally isomorphic to 
$\big(\mathrm{PSL}(2,\mathbb{R})\big)^k$ with $k\le n$. Note that $\rho^*$
factors through
$$
 H^{2n}_c(\mathrm{GL}^+(2n,\mathbb{R}))
 \longrightarrow H^{2n}_c(S)\cong H^{2n}_c\left(\Pi_{i=1}^k 
 \mathrm{PSL}(2,\mathbb{R})\right)\longrightarrow
 H^{2n}(\Gamma).
$$ 
If $k<n$, then the middle
cohomology group is zero and hence $\rho^*(\varepsilon_{2n})$
vanishes. If $k=n$, then Proposition \ref{prop:rigid} gives us
further that either $t(S)>0$, in which case $\rho^*(\varepsilon_{2n})=0$ (Lemma \ref{lem:trho=0}), or
$S$ is conjugate to a subgroup of $\text{Bl}_{2,n}^+$.

\medskip

\noindent{\it Case 3: The mixed case.}


Applying Proposition \ref{prop:mixed} and arguing as above, we derive that if either $(i),(ii)$ or $(iii')$ of \ref{prop:mixed} holds then $\rho^*(\varepsilon_{2n})=0$. Thus we may assume in this case as well that $\rho(\gC)$ is conjugate to a subgroup of $\text{Bl}_{2,n}^+$.

\medskip

We have therefore reduced to the situation where, up to conjugating the image, $\rho$ factors through a map
$\rho_0$ into $\text{Bl}_{2,n}^+$:
\[
\xymatrix{\Gamma\ar[r]^{\rho \ \ \ \ \ }\ar[dr]_{\rho_0} & \GL^+(2n,\bbr) \\
                 & \text{Bl}_{2,n}^+.\ar@{^{(}->}[u]^i      }
\]
The Whitney product formula for the Euler class gives
$$
 i^*(\varepsilon_{2n})=\varepsilon_2\cup\ldots\cup \varepsilon_2 \in
 H^{2n}_c(\text{Bl}_{2,n}^+),
$$
and hence
$$
 \left\Vert \rho^{\ast}(\varepsilon_{2n})\right\Vert _{\infty}=\left\Vert \rho_0^{\ast}\circ i^*(\varepsilon_{2n})\right
 \Vert_{\infty} =\left\Vert \rho_0^{\ast}(\varepsilon_{2}\cup \ldots \cup \varepsilon_2)\right\Vert _{\infty} \leq\left\Vert
 \varepsilon_{2}\cup\ldots\cup\varepsilon_{2}\right\Vert _{\infty}.
$$
\end{proof}


\section{The proofs of Theorems \ref{Theorem: Milnor-Wood for products}, \ref{Theorem: sharp}, \ref{thm: gnrl Euler number}}\label{section: Proof of the main theorems}

\noindent{\bf Proof of Theorem \ref{Theorem: Milnor-Wood for products}.} 
Let $M$ be a\ closed oriented Riemannian manifold with universal cover
$\mathcal{H}^n=\prod_{i=1}^n\mathcal{H}$. Let $\Gamma$ be the fundamental
group of $M$ embedded as a cocompact lattice in
$\text{Isom}(\mathcal{H}^n)^+$ acting on $\mathcal{H}^n$ by deck transformations.
Let $\xi$ be a $\mathrm{GL}^+(2n,\mathbb{R}%
)$-bundle over $M$. 
Suppose that $\xi$ admits a flat structure, and
let $\rho:\Gamma\rightarrow\mathrm{GL}^+(2n,\mathbb{R})$ be the corresponding
representation. 
Identifying $H^{2n}(M)$ with $H^{2n}(\Gamma)$, $\varepsilon(\xi)$ considered as an element of $H^{2n}(\Gamma)$ is equal to $\rho^{\ast}(\varepsilon_{2n})$.
Since this identification is an isometry,
Theorem \ref{Thm: norm of the Euler class} gives
\[
\left\Vert \varepsilon_{2n}(\xi)\right\Vert _{\infty}=  \left\Vert \rho^*(\varepsilon_{2n})\right\Vert _{\infty}  \leq\left\Vert \varepsilon
_{2}\cup\ldots\cup\varepsilon_{2}\right\Vert _{\infty}.
\]
Combining Equation (\ref{Equ: b(M)=product of norms}), the proportionality principle
established in Proposition \ref{prop: Proportionality Principle
Simpl Vol and Euler} and the last inequality, we conclude
\[
\left\vert \left\langle \varepsilon(\xi),[M]\right\rangle \right\vert
=\left\Vert \varepsilon(\xi)\right\Vert _{\infty}\left\Vert M\right\Vert
=\frac{\left\Vert \varepsilon(\xi)\right\Vert _{\infty}}{\left\Vert
\varepsilon_{2}\cup\ldots\cup\varepsilon_{2}\right\Vert _{\infty}}\frac{\chi
(M)}{(-2)^{n}}\leq\frac{1}{\left(  -2\right)  ^{n}}\chi(M),
\]
which completes the proof of the theorem.
\qed

\medskip

\noindent{\bf Proof of Theorem \ref{Theorem: sharp}.}
Let $\Gamma<\text{Isom}(\mathcal{H}^n)^+$ be the fundamental group of
$M$ with $\ti\gC$ a cofaithful lift of $\gC$
in $G_n^+$
and $\rho:\gC\to\ti\gC< G_n^+$ the cofaitfhul map. Up to replacing $M$ by a finite cover, which amounts to multiplying each side of the equality we are proving by the same number (the degree of the cover), we may assume that $\gC\le\prod_{i=1}^n\PSL(2,\BR)$ and $\ti\gC\le\prod_{i=1}^n\SL(2,\BR)$.
We thus have a commutative diagram
\[
\xymatrix{\Gamma\ar@{^{(}->}[r] \ar@{^{(}->}^{\rho \ \ \ \ \ \ \ }[rd] & \Pi_{i=1}^n\PSL(2,\bbr) \\
          & \Pi_{i=1}^n\SL(2,\bbr).\ar@{->>}[u]      }
\]
Since $\rho$ is injective and $\rho(\gC)$ is a cocompact lattice
in $\Pi_{i=1}^n\SL(2,\bbr)$, the
induced map $\rho^*:H^*_c(\Pi_{i=1}^n\SL(2,\bbr))\rightarrow
H^*(\gC)$ is an isometric embedding (cf. \cite[Proposition
8.6.2]{Mo01} or \cite[Theorem 3]{Bu06}). 
Thus
$$
 \| \rho^*(\varepsilon_2\cup \ldots\cup \varepsilon_2)\|_\infty=\| \varepsilon_2\cup \ldots
 \cup \varepsilon_2 \|_\infty.
$$ 
Let $\xi_\rho$ be the flat
$\GL^+(2n,\bbr)$-bundle over $M$ corresponding to the representation
$\rho:\gC\rightarrow \Pi_{i=1}^n \SL(2,\bbr) 
\cong \text{SBl}_{2,n}(\BR)<\GL^+(2n,\bbr)$. Then
$$
 |\chi(\xi_\rho)|=|\langle \rho^*(\varepsilon_2\cup
 \ldots\cup \varepsilon_2),[M]\rangle |= \| \varepsilon_2\cup
 \ldots\cup\varepsilon_2\|_\infty \|M\|=\frac{1}{2^n}|\chi(M)|,$$ where the
last equality follows from Proposition \ref{prop: Proportionality
Principle Simpl Vol and Euler}.
\qed

\medskip

\noindent{\bf Proof of Theorem \ref{thm: gnrl Euler number}.}
Let $N$ be a closed $\mathcal{H}^n$-manifold of the form
$$N=\Sigma_{g_1}\times\ldots\times \Sigma_{g_k}\times N',$$
where $N'$ is rigid, $k\ge 0$, and the
$\Sigma_{g_i}$'s are surfaces of genus $g_i\geq 2$, and let $\xi_N$ be a flat  $\GL^+(2n,\mathbb{R})$-bundle over $N$. In order to show that 
$$
\chi(\xi_N)\in \bigg\{ \pm \frac{\chi(N')}{2^{\mathrm{Dim}(N')/2}}\Pi_{i=1}^k \ell_i :  |\ell_i|\leq g_i-1\}\cup\{0\bigg\},
$$
we split the proof into three different cases, as before:

\medskip
\noindent{\it Case 1: $N=\Sigma_{g_1}\times\cdots\times \Sigma_{g_n}$ is a product of surfaces and $N^\prime$ is trivial.} The fundamental group $\Gamma$ of
$N$ is a product of $n$ surface groups $\gC_i$. Let $\rho:\Gamma\rightarrow\GL^+(2n,\bbr)$ be the oriented
representation corresponding to $\xi_N$. If the restriction of $\rho$ to any of the surface
groups is amenable, then $\rho^*(\varepsilon_{2n})=0$ by Lemma
\ref{lemma: rho(Gamma) amenable}, and if $t(\rho)>0$ then again $\rho^*(\varepsilon_{2n})=0$ by Lemma \ref{lem:trho=0}, so let us assume that this is not
the case. 
It now follows from Proposition \ref{prop: Case 1 Gamma completely reducible},
that $\rho (\gC)$ is conjugate to a subgroup of $\mathrm{Bl}^+_{2,n}$. Since
conjugation induces isomorphisms of cohomology groups, we can without loss
of generality assume that the image of $\rho$ is in fact contained
in $\mathrm{Bl}^+_{2,n}$. Moreover, it follows from Proposition \ref{prop: Case 1 Gamma completely reducible} that up to reordering the factors if necessary, the $i$-th factor $\Gamma_i$ is mapped irreducibly to the $i$-th factor of $\mathrm{Bl}^+_{2,n}\cong \Pi_{i=1}^n\GL^+(2,\mathbb{R})$ and scalarly to every other factor. Denote by $r_i:\Gamma_i\rightarrow \GL^+(2,\mathbb{R})$ the restriction of $\rho$ to $\Gamma_i$ composed with the projection onto the $i$-th factor of $\mathrm{Bl}^+_{2,n}$. 
Then $\rho$ and $(r_1,\ldots,r_n)$ have the same projection to $ \Pi_{i=1}^n\PSL(2,\mathbb{R})$. Thus by Lemma \ref{lem: reps with same proj to products of PSL(2,R)},
$$\rho^*(\varepsilon_{2n})=r_1^*(\varepsilon_2)\cup\ldots\cup r_n^*(\varepsilon_2).$$
We hence obtain that 
\begin{eqnarray*}
\chi(\xi_N)=\langle \varepsilon_{2n}, [N] \rangle &=& \langle r_1^*(\varepsilon_2)\cup\ldots\cup r_n^*(\varepsilon_2), [\Sigma_{g_1}\times \ldots\times \Sigma_{g_n}] \rangle \\
&=& \langle r_1^*(\varepsilon_2),[\Sigma_{g_1}]  \rangle \cdot \ldots \cdot \langle r_n^*(\varepsilon_2), [\Sigma_{g_n}] \rangle.
\end{eqnarray*}
Finally, for every $i$, $\langle r_i^*(\varepsilon_2), [\Sigma_{g_i}] \rangle$ is the Euler number of a flat bundle over the surface $\Sigma_{g_i}$, and hence, by Equation (\ref{Equ: b(M)=product of norms}) and Proposition \ref{prop: Proportionality Principle Simpl Vol and Euler}, satisfies
$$
 |\langle r_i^*(\varepsilon_2), [\Sigma_{g_i}] \rangle|=\| r_i^*(\varepsilon_2)\|_\infty\|\Sigma_{g_i}\|\le \frac{\| r_i^*(\varepsilon_2)\|_\infty}{\|\varepsilon_2\|_\infty}\frac{\chi(\gS_{g_i})}{-2}\le g_i-1,
$$
which is Milnor's classical 
inequality \cite{Mi58}.

\medskip

\noindent{\it Case 2: $N=N^\prime$ is rigid.}
Let $\Gamma$ denote the fundamental group of $N$ and let
$\rho:\Gamma\rightarrow \GL^+(2n,\bbr)$ be a representation inducing
the flat bundle $\xi_N$. Proposition \ref{prop: Proportionality
Principle Simpl Vol and Euler} gives
$$|\chi(\xi_N)|=\|\rho^*(\varepsilon_{2n}) \|_{\infty} \|M\|=\frac{\|\rho^*(\varepsilon_{2n}) \|_{\infty}}{\|\varepsilon_{2}\cup \ldots \cup\varepsilon_{2} \|_{\infty}}\frac{1}{2^n}|\chi(M)|.$$
Hence the result will follow from the next claim:

\medskip
\noindent {\bf Claim.} Either $\rho^*(\varepsilon_{2n})=0$ or $\|\rho^*(\varepsilon_{2n}) \|_{\infty}=\|\varepsilon_{2}\cup \ldots \cup\varepsilon_{2} \|_{\infty}$.
\medskip

Since the inclusion of a finite index subgroup in $\Gamma$ induces isometric embeddings on the 
standard and 
bounded cohomology groups, we can without loss of generality replace $\Gamma$ by a finite index subgroup. 
Thus, as in the proof of Theorem \ref{Thm: norm of the Euler class}, applying Proposition \ref{prop:rigid}, we can reduce to the situation where either $\rho^*(\varepsilon_{2n})=0$, or $\rho$ is injective, and its image $\rho(\Gamma)$ is, up to conjugation, contained, discrete and cocompact in
$S=\text{SBl}_{2,n}$. In the latter case,  we again invoke \cite[Proposition
8.6.2]{Mo01} (or \cite[Theorem 3]{Bu06}) to conclude that the
induced map $\rho^*:H^*_c(\text{SBl}_{2,n})\rightarrow H^*(\Gamma)$
is an isometric embedding and
$$\|\rho^*(\varepsilon_{2n}) \|_{\infty}=\|\varepsilon_{2}\cup\ldots\cup\varepsilon_{2} \|_{\infty},  $$
which finishes the proof of the claim. 

\medskip

\noindent{\it Case 3: The mixed case:} $N=\Sigma_{g_1}\times \ldots\times \Sigma_{g_k}\times N'$ with $k>0$ and $N^{\prime}$ nontrivial. The fundamental group $\Gamma$ of $N$ decomposes as $\Gamma=\Gamma_1 \times \ldots \times \Gamma_k \times \Gamma^\prime$, where the $\Gamma_i$'s are surface groups and $\Gamma^\prime$ is the fundamental group of the rigid factor $N^\prime$. Let
$\rho:\Gamma\rightarrow \GL^+(2n,\bbr)$ be a representation inducing
the flat bundle $\xi_N$. 


Replacing the rigid factor $\Gamma^\prime$ by a finite index subgroup has no effect on the equality we are proving (both sides are multiplied by the index of the subgroup). Hence, assuming that $\rho^*(\varepsilon_n)\ne 0$ , arguing as above, we may reduce the situation to case $(iii')$ of Proposition \ref{prop:mixed}.


Furthermore, by Proposition \ref{prop:prod}, up to conjugation we can assume that the representation  
$$
 \rho:\Gamma=\Gamma_1 \times \ldots \times \Gamma_k \times \Gamma^\prime \rightarrow \mathrm{Bl}^+_{2,n}\cong   \GL^+(2,\BR)\times \ldots \times \GL^+(2,\BR)
$$ 
maps surface group factor $\Gamma_i$ irreducibly to the $i$-th factor of $\mathrm{Bl}^+_{2,n}$ and scalarly to the others and, the rigid factor $\Gamma^\prime$ is mapped scalarly to each of the $k$ first factors of $\mathrm{Bl}^+_{2,n}$. Define $r_i:\Gamma_i\rightarrow \GL^+(2,\BR)$ as the restriction of $\rho$ to $\Gamma_i$ composed with the projection onto the $i$-th factor of $\mathrm{Bl}^+_{2,n}$ and $r^\prime:\Gamma^\prime \rightarrow \Pi_{i=k+1}^n \GL^+(2,\BR)$ as the restriction of $\rho$ to $\Gamma^\prime$ composed with the projection onto the last $n-k$ factors of $\mathrm{Bl}^+_{2,n}$. 
Then $\rho$ and 
$$
 (r_1,\ldots,r_k,r^\prime):\Gamma_1 \times \ldots\times \Gamma_k\times \Gamma^\prime\rightarrow \mathrm{Bl}^+_{2,n}
$$ 
have the same projection to  $\Pi_{i=1}^n \PSL(2,\BR)$, hence by Lemma \ref{lem: reps with same proj to products of PSL(2,R)} the corresponding pullbacks of the Euler class coincide.
Thus the first part of the theorem follows in this case from its validity in cases 1 and 2. 

\medskip
Let us now explain the second statement of the theorem, namely that if $N$ is cofaithful, then all the integers satisfying the given rule are actually attained as Euler numbers of flat bundles. The assertion is true for surfaces by \cite{Mi58}. It is further true for any rigid cofaithful manifold $N^\prime$, since the Euler number is then equal to either $0$ or $\pm \chi(N^\prime)/2^{\dim(N^\prime)/2}$, where $0$ is realized for example by the trivial bundle, and the latter, up to orientation, by the cofaithful map (Theorem \ref{Theorem: sharp}). The general case follows by considering direct products of flat bundles over the factors of $N$ with the appropriate Euler numbers.
\qed





\section{The proof of Theorem \ref{thm:rigidity}}\label{section: pf of thm for rigid mflds}
For a group $G$ we denote by $S_n\ltimes G^n$ the semidirect product where $S_n$ permutes the factors of $G^n=\prod_{i=1}^nG$. When $G$ is a subgroup of $\GL(2,\BR)$, $S_n\ltimes G^n$ admits a natural embedding in $\GL(2n,\BR)$ where $G^n$ embeds diagonally in $\text{Bl}_{2,n}$. We will not distinguish between $S_n\ltimes G^n$ and this representation of it, and denote by $(S_n\ltimes G^n)^+$ its intersection with $\GL^+(2n,\BR)$. In particular $G_n^+=(S_n\ltimes\GL^1(2,\BR))^+$, the $2^n$ sheeted cover of $\text{Isom}(\mathcal{H}^n)^+$, is embedded in this form.

Recall that $s$ denotes the quotient map 
$$
 s:\text{Isom}(\mathcal{H}^n)\cong S_n\ltimes\PGL(2,\BR)^n\to S_n. 
$$
For a subgroup $\gC\le \text{Isom}(\mathcal{H}^n)$ and an abelian group $A$ we denote by $H^1_\text{Sym}(\gC,A^n)$ the first cohomology group of $\gC$ in $A^n$ with respect to the action of $\gC$ on $A^n$ that permutes the factors according to the map $s$.
 

\begin{lem}\label{lem:H^1-finite}
Suppose that $\gC\le \text{Isom}(\mathcal{H}^n)$ is a rigid lattice. Then: 
\begin{enumerate}
\item $ H^1_\text{Sym}(\gC,\BR^n)=\langle 0\rangle$.
\item The inclusion of $\langle\pm 1\rangle$ in $\BR^*$ induces an isomorphism $H^1_\text{Sym}(\gC,(\langle\pm 1\rangle)^n)\cong H^1_\text{Sym}(\gC,(\BR^*)^n)$, and this cohomology group is finite.
\end{enumerate}
\end{lem}

\begin{proof}
Let $\gC_0$ be the intersection of $\gC$ with the identity connected component of $\text{Isom}(\mathcal{H}^n)$. Then $\gC_0$ is a rigid lattice in $\PSL(2,\BR)^n$, hence has finite abelianization and in particular, no nontrivial homomorphisms to $\BR$, as follows for instance by Margulis' normal subgroups theorem. Moreover, as $s|_{\gC_0}$ is trivial, any element of $H^1_\text{Sym}(\gC,\BR^n)$ restricts to a homomorphism on $\Gamma_0$, and is hence identically $0$ on $\gC_0$. It follows from the cocycle equation that any element $\ga\in H^1_\text{Sym}(\gC,\BR^n)$ is constant on cosets of $\gC_0$:
$$
 \ga(\gc_0\gc)=\ga(\gc_0)+s(\gc_0)\cdot\ga(\gc)=\ga(\gc),~\forall \gc_0\in\gC,\gc\in\gC.
$$
Thus $H^1_\text{Sym}(\gC,\BR^n)\cong H^1_\text{Sym}(\gC/\gC_0,\BR^n)$, which is trivial since $\gC/\gC_0$ is finite and hence has property $(T)$. This proves $(1)$.

Since $\BR^*\cong \BR\oplus\BZ/(2)$ we have 
$$
 H^1_\text{Sym}(\gC,(\BR^*)^n)\cong H_\text{Sym}^1(\gC,\BR^n)\oplus H^1_\text{Sym}(\gC,(\BZ/(2))^n),
$$
and the first statement of $(2)$ follows from $(1)$.

Moreover, since $\gC_0$ is a lattice in $\prod_{i=1}^n\PSL(2,\BR)$, it is finitely generated, and hence $\text{Hom}(\gC_0,(\BZ/(2))^n)$ is finite. By the cocycle equation, an element $\ga\in H^1_\text{Sym}(\gC,(\BZ/(2))^n)$ is determined by its restriction to $\gC_0$ and by its values on a finite set of cosets representatives for $\gC_0$ in $\gC$. 
Thus $H^1_\text{Sym}(\gC,(\BZ/2)^n)$ is finite. 
\end{proof}

We now give the proof of Theorem \ref{thm:rigidity}: Let $M$ be a closed locally $\mathcal{H}^n$ rigid manifold with fundamental group $\gC$.

\medskip

\noindent {\bf Nonzero Euler number implies cofaithfulness:}
Suppose that $M$ admits a $\dim(M)$-dimensional flat vector bundle with nonzero Euler number. Up to reversing the orientation we may assume that the Euler number is positive, hence by Corollary \ref{cor: Euler number for rigid}, equals $|\frac{1}{2^n}\chi(M)|$. Let $\rho:\gC\to\GL^+(2n,\BR)$ be the linear representation inducing this structure. By Proposition \ref{prop:rigid}, $\gC$ admits a normal finite index subgroup $\gD$ such that $\rho(\gD)$ is contained in a semisimple Lie group $S$ with at most $n$ simple factors, for which every noncompact simple factor is locally isomorphic to $\PSL(2,\BR)$. 
If $S$ admits less than $n$ factors or a compact factor,  Lemma \ref{lemma: can mod out amenable} combined with the fact that $S_{nc}$ (the product of noncompact factors of $S$) is locally isomorphic in that case to the isometry group of a symmetric space of dimension strictly smaller than $2n$, yields that $\rho^*(\varepsilon_{2n})=0$, contrary to our assumption. Therefore, $S$ is locally isomorphic to $\PSL(2,\BR)^n$.
By Lemma \ref{lem:trho=0}, $t(S)=0$, hence we may assume that $(b_1),(b_2),(b_3)$ of Proposition \ref{prop:rigid} hold for $\gD$ and $\rho$, and in particular, up to replacing $\rho$ by some conjugate representation, that $S=\text{SBl}_{2,n}$ and $\rho(\gD)$ is a cocompact lattice there. 
In particular, since $\rho(\gD)$ is normal in $\rho(\gC)$ and, by $(b3)$ and Borel's density theorem, Zariski dense in $S$, it follows that $\rho(\gC)$ is contained in the normalizer of $S=\text{SBl}_{2,n}$, namely, in $(S_n\ltimes\GL(2,\BR)^n)^+$. By replacing $\gD$ further by some characteristic finite index subgroup if necessary, we may assume additionally that $\gD$ itself lies in $\PSL(2,\BR)^n$, the identity connected component of $\text{Isom}^+(\mathcal{H}^n)\cong (S_n\ltimes \PGL(2,\BR)^n)^+$.

Let $f:\gD\to \prod_{i=1}^n\PSL_2(\BR)$ be the composition of $\rho|_\gD$ followed by $\Ad:\text{SBl}_{2,n}\to\Ad(\text{SBl}_{2,n})\cong \prod_{i=1}^n\PSL_2(\BR)$. By Mostow's rigidity theorem $f$ extends to an isomorphism $\ti f: \prod_{i=1}^n\PSL_2(\BR)\to \prod_{i=1}^n\PSL_2(\BR)$ which, by reordering the factors of the target group, preserves the order of the factors. For every $\gc\in\gC$ and $\gd\in\gD$ we have
$$
 f(\gc\gd\gc^{-1})=\rho(\gc)f(\gd)\rho(\gc)^{-1}
$$
where the term on the right hand side is understood as $\Ad (\rho(\gc)\rho(\gd)\rho(\gc)^{-1})$. 
Hence by Borel's density theorem we have 
\begin{equation}\label{eq:22222}
 \ti f (\gc g\gc^{-1})=\rho(\gc)\ti f(g)\rho(\gc)^{-1}
\end{equation}
for any $g\in\PSL(2,\BR)^n$. In particular, $s(\rho(\gc))=s(\gc)$, i.e. $\rho(\gc)$ and $\gc$ induce the same permutation on the factors of $\SL(2,\BR)^n$ and $\PSL(2,\BR)^n$ respectively. Thus we can write $\rho(\gc)\in (S_n\ltimes\GL(2,\BR)^n)^+$ as 
$$
 (s(\gc),\rho_1(\gc),\ldots,\rho_n(\gc)),
$$
where $\rho_i(\gc)$ is the component of $\rho(\gc)$ in the $i$'th factor of $\GL(2,\BR)^n$.
One can easily verify that the map $\ga:\gC\to\BR^n$ given by
$$
 \ga(\gc):=(\log\det (\rho_1(\gc)^2),\ldots,\log\det (\rho_n(\gc)^2))
$$
is a cocycle, and hence by Lemma \ref{lem:H^1-finite} cohomologous to the trivial map. Thus there exists a diagonal matrix of the form
$\gL=\text{diag}(\gl_1,\gl_1,\ldots,\gl_n,\gl_n)$ such that the image of $\rho$ conjugated by $\gL$ is contained in $G_n^+=(S_n\ltimes \GL^1(2,\BR)^n)^+$.

Finally let $\psi:(S_n\ltimes \GL^1(2,\BR)^n)^+\to (S_n\ltimes \PGL(2,\BR)^n)^+$ be the canonical projection determined by modding out the center of $\GL^1(2,\BR)^n$. Since $(S_n\ltimes \PGL(2,\BR)^n)^+$ is isomorphic to $\Aut^+(\PSL(2,\BR)^n)$, where the faithful action of the first on $\PSL(2,\BR)^n$ is by conjugation on its identity connected component, we deduce from Equation \ref{eq:22222} that $\psi\circ\rho:\gC\to(S_n\ltimes \PGL(2,\BR)^n)^+$ is injective.
Hence 
$$
 \rho\circ (\psi\circ\rho)^{-1}:\psi(\rho(\gC))\to (S_n\ltimes \GL^1(2,\BR)^n)^+
$$
induces a faithful lift, showing that $M$ is cofaithful.


\vspace{1cm}

\noindent{\bf Characterisation of flat structures with positive Euler number as elements in $H^1_\text{Sym}(\gC,\BZ/(2)^n)$:}
Suppose now that $M$ is cofaithful, and let $\rho_1:\gC\to G_n^+\le\GL^+(2n,\BR)$ be the representation induced by the cofaithful lift of $\gC$ from $S_n\ltimes(\PGL(2,\BR)^n)^+$ to $(S_n\ltimes (\GL^1_2)^n)^+\cong G_n^+$ where the latter is realized as the subgroup $\GL(2n,\BR)^+$ described at the beginning of this section.
Let $\rho:\gC\to\GL^+(2n,\BR)$ be another representation, and suppose that $\rho^*(\varepsilon_{2n})\ne 0$. 
As above, we derive from Proposition \ref{prop:rigid} that there is a finite index normal subgroup $\gD\le\gC$ such that the restriction of $\rho$ to $\gD$ is faithful and $\rho(\gD)$ is conjugated to some cocompact lattice in $\text{SBL}_{2,n}(\BR)$. Thus, it follows from Mostow's rigidity theorem that $\rho^g|_\gD=\rho_1|_\gD$ for an appropriate element $g\in \GL(2n,\BR)$ (by $\rho^g$ we mean $\rho$ composed with the conjugation by $g$). Now let $\gc$ be an arbitrary element of $\gC$. Since $\gD$ is normal in $\gC$, we have 
\begin{eqnarray*}
 \rho_1(\gc^{-1})\rho_1(\gd)\rho_1(\gc)&=&\rho_1(\gc^{-1}\gd\gc)=\rho^g(\gc^{-1}\gd\gc)=\rho^g(\gc^{-1})\rho^g(\gd)\rho^g(\gc)\\
 &=&\rho^g(\gc^{-1})\rho_1(\gd)\rho^g(\gc),
\end{eqnarray*}
for every $\gd\in\gD$. Thus $\rho^g(\gc)\rho_1(\gc)^{-1}$ lies in the centralizer of $\rho_1(\gD)$, and by
Borel's density theorem, in the centralizer of $\text{SBL}_{2,n}(\BR)$.
Hence $\rho^g(\gc)$ is equal to $\rho_1(\gc)$ multiplied from the left by a matrix of the form 
$$
 \text{diag}(\chi_1(\gc),\chi_1(\gc),\ldots,\chi_n(\gc),\chi_n(\gc))
$$ 
where the $\chi_1,\ldots,\chi_n$ are functions on $\gC$ taking values in $\BR^*$. It also follows that the permutation representation determined by conjugating by $\rho^g(\gc)$ on the factors of $\text{SBl}_{2,n}$ is the same as the one coming from conjugation by $\rho_1(\gc)$, namely that $s(\rho(\gc))=s(\gc)$.
A simple calculation shows that the map $\ga:\gC\to(\BR^*)^n$ given by
$\ga(\gc)=(\chi_1(\gc),\ldots,\chi_n(\gc))$ is a cocycle. Hence by Lemma \ref{lem:H^1-finite}, after conjugating $\rho^g$ further by some diagonal matrix of the form 
$\text{diag}(\gl_1,\gl_1,\ldots,\gl_n,\gl_n)$ if necessary, we may assume that $\ga$ takes its values in $\langle\pm 1\rangle^n$. It follows that the conjugacy class of $\rho$ in $\GL(2n,\BR)^+$ is completely determined by the cocycle $\ga$ and the sign of the determinant of $g$, where $\det (g)<0$ corresponds to negative Euler number.  

Let $\gC_0$ be as in the proof of Lemma \ref{lem:H^1-finite}. Then the restriction of every $[\ga]\in H^1_s(\gC,\langle \pm 1\rangle^n)$ to $\gC_0$ is a homomorphism. Set
$$
 \gC^0=\gC_0\cap \{\ker [\ga]:[\ga]\in H^1_s(\gC,\langle \pm 1\rangle^n).
$$
Then $\gC^0$ is a characteristic subgroup of finite index in $\gC$ and the restrictions to $\gC^0$ of all the representations of $\gC$ which induce flat vector bundles with nonzero (resp. positive) Euler number are conjugate in $\GL(2n,\BR)$ (resp. in $\GL^+(2n,\BR)$). 

For the other direction, any cocycle $\ga:\gC\to\langle\pm 1\rangle^n$ 
produces such a representation by setting
$$
 \rho_\ga(\gc):=\text{diag}\big(\ga(\gc)_1,\ga(\gc)_1,\ldots,\ga(\gc)_n,\ga(\gc)_n\big)\rho_1(\gc)
$$
which gives rise to a flat vector bundle with positive Euler number.
Moreover if $\ga$ and $\gb$ are two such cocycles and the corresponding representations $\rho_\ga,\rho_\gb$ are conjugate in $\GL^+(2n,\BR)$, the conjugating element must centralize $\rho_\ga(\gC^0)$ and hence, by the Borel density theorem, belongs to the centralizer of $\text{SBl}_{2,n}(\BR)$, i.e. it is of the form $\text{diag}(\gl_1,\gl_1,\ldots,\gl_n,\gl_n)$. This shows that $\ga$ and $\gb$ are cohomologus as cocycles with coefficients in $(\BR^*)^n$, namely, they represent the same element in $H^1_\text{Sym}(\gC,(\BR^*)^n)$. By Lemma \ref{lem:H^1-finite} (2), they represent the same element in $H^1(\gC,\langle\pm 1\rangle)$. 
This completes the proof of the theorem.
\qed

\begin{rem}
The finiteness of the number of flat vector bundles of any given dimension over any rigid manifold is true in general. Indeed, one can deduce from the proof of Margulis superrigidity theorem, the fact that any rigid lattice has finite abelianization, and Jordan's theorem, that $\pi_1(M)$ admits only finitely many nonequivalent representations in any fixed dimension. 
\end{rem}


\noindent{\textrm{Acknowledgments}:} \emph{We thank Bill Goldman for pointing out some essential references to us. We thank Etienne Ghys for many useful comments and especially for showing us Example \ref{Example Ghys}. 
We thank Uri Bader for several helpful remarks and in particular for showing us an elegant proof of Lemma \ref{lem:H^1-finite}. 
Finally, we thank Peter Storm for his interest in this work, for hours of conversations and for his helpful comments which definitely improved this paper.}


\end{document}